\documentclass[notitlepage,11pt,reqno]{amsart}
\usepackage[foot]{amsaddr}
\usepackage{amssymb,nicefrac,bm,upgreek,mathtools,verbatim,enumerate}
\usepackage[final]{hyperref}
\usepackage[mathscr]{eucal}
\usepackage{dsfont}
\usepackage[normalem]{ulem}
\usepackage{amsopn}
\usepackage{color}
\usepackage[margin=1in]{geometry}
\newcommand{\ttup}[1]{\textup{(}#1\textup{)}}
\newcommand{\stkout}[1]{\ifmmode\text{\sout{\ensuremath{#1}}}\else\sout{#1}\fi}
\newtheorem{lemma}{Lemma}[section]
\newtheorem{theorem}{Theorem}[section]

\theoremstyle{definition}
\newtheorem{definition}{Definition}[section]
\newtheorem{algorithm}{Algorithm}[section]
\newtheorem{assumption}{Assumption}[section]

\newtheorem{example}{Example}[section]

\theoremstyle{remark}
\newtheorem{remark}{Remark}[section]
\numberwithin{theorem}{section}
\numberwithin{equation}{section}

\hypersetup{
  colorlinks=true,
  citecolor=blue,
  linkcolor=blue,
  frenchlinks=false,
  breaklinks}
\usepackage[abbrev,msc-links,nobysame]{amsrefs} 
\usepackage[capitalize,nameinlink]{cleveref}

\crefname{section}{Section}{Sections}
\crefname{subsection}{Section}{Sections}
\crefname{condition}{Condition}{Conditions}
\crefname{hypothesis}{Hypothesis}{Conditions}
\crefname{assumption}{Assumption}{Assumptions}
\crefname{lemma}{Lemma}{Lemmas}
\crefname{fact}{Fact}{Facts}

\Crefname{figure}{Figure}{Figures}

\crefformat{equation}{\textup{#2(#1)#3}}
\crefrangeformat{equation}{\textup{#3(#1)#4--#5(#2)#6}}
\crefmultiformat{equation}{\textup{#2(#1)#3}}{ and \textup{#2(#1)#3}}
{, \textup{#2(#1)#3}}{, and \textup{#2(#1)#3}}
\crefrangemultiformat{equation}{\textup{#3(#1)#4--#5(#2)#6}}%
{ and \textup{#3(#1)#4--#5(#2)#6}}{, \textup{#3(#1)#4--#5(#2)#6}}%
{, and \textup{#3(#1)#4--#5(#2)#6}}

\Crefformat{equation}{#2Equation~\textup{(#1)}#3}
\Crefrangeformat{equation}{Equations~\textup{#3(#1)#4--#5(#2)#6}}
\Crefmultiformat{equation}{Equations~\textup{#2(#1)#3}}{ and \textup{#2(#1)#3}}
{, \textup{#2(#1)#3}}{, and \textup{#2(#1)#3}}
\Crefrangemultiformat{equation}{Equations~\textup{#3(#1)#4--#5(#2)#6}}%
{ and \textup{#3(#1)#4--#5(#2)#6}}{, \textup{#3(#1)#4--#5(#2)#6}}%
{, and \textup{#3(#1)#4--#5(#2)#6}}

\crefdefaultlabelformat{#2\textup{#1}#3}
\newcommand{\vertiii}[1]{{\left\vert\kern-0.25ex\left\vert\kern-0.25ex\left\vert #1 
    \right\vert\kern-0.25ex\right\vert\kern-0.25ex\right\vert}}
\newcommand{\lamstr}{\lambda^{\mspace{-2mu}*}}% lambda with corrected 'star'
\newcommand{\process}[1]{{\{#1_t\}_{t\ge0}}}

\newcommand{\Uadm}{\mathfrak U}
\newcommand{\Act}{\mathbb{U}}
\newcommand{\Usm}{\mathfrak U_{\mathsf{sm}}}
\newcommand{\bUsm}{\overline{\mathfrak U}_{\mathsf{sm}}}
\newcommand{\Usms}{\mathfrak U^*_{\mathsf{sm}}}
\newcommand{\cA}{{\mathcal{A}}}  % generator
\newcommand{\sB}{{\mathscr{B}}}  % Ball
\newcommand{\cC}{{C}}   % Continuous functions
\newcommand{\sE}{{\mathscr{E}}} 
\newcommand{\sF}{{\mathfrak{F}}}   % sigma field
\newcommand{\cG}{{\mathcal{G}}}  % operator
\newcommand{\cK}{{\mathcal{K}}}  % compact set
\newcommand{\sK}{{\mathscr{K}}}  % 
\newcommand{\sL}{{\mathscr{L}}}  % Linear operator
\newcommand{\Lp}{{L}}            % Lp
\newcommand{\Lpl}{L_{\text{loc}}}            % Lp local
\newcommand{\cM}{{\mathcal{M}}}  % compact set

\newcommand{\Lyap}{{\mathscr{V}}}  % Lyapunov
\newcommand{\cX}{{\mathcal{X}}}

\newcommand{\RR}{\mathds{R}}
\newcommand{\NN}{\mathds{N}}

\newcommand{\Rd}{{\mathds{R}^{d}}}
\DeclareMathOperator{\Exp}{\mathbb{E}}
\DeclareMathOperator{\Prob}{\mathbb{P}}
\newcommand{\D}{\mathrm{d}}
\newcommand{\E}{\mathrm{e}}
\newcommand{\Ind}{\mathds{1}}   % indicator function

\newcommand{\Sob}{{\mathscr W}}    % Sobolev Space
\newcommand{\Sobl}{{\mathscr W}_{\mathrm{loc}}} % Sobolev Space(local)

\newcommand{\df}{\coloneqq}
\newcommand{\transp}{^{\mathsf{T}}}

\DeclareMathOperator*{\trace}{Tr}

\DeclareMathOperator*{\Argmin}{Arg\,min}
\DeclareMathOperator*{\Argmax}{Arg\,max}

\newcommand{\sorder}{{\mathfrak{o}}}
\newcommand{\grad}{\nabla}
\newcommand{\uuptau}{{\Breve\uptau}}
\newcommand{\sEmin}{{\mathscr{E}^*_{\mathrm{min}}}} 
\newcommand{\sEmax}{{\mathscr{E}^*_{\mathrm{max}}}} 

\newcommand{\tv}{{\rule[-.45\baselineskip]{0pt}{\baselineskip}\mathsf{TV}}}

\newcommand{\abs}[1]{\lvert#1\rvert}
\newcommand{\norm}[1]{\lVert#1\rVert}
\newcommand{\babs}[1]{\bigl\lvert#1\bigr\rvert}

\newcommand{\bnorm}[1]{\bigl\lVert#1\bigr\rVert}

\begin{document}
\title[On the policy improvement algorithm for risk-sensitive control]
{On the policy improvement algorithm\\[5pt] for ergodic risk-sensitive control}

\author[Ari Arapostathis]{Ari Arapostathis$^{\dag}$}
\address{$^{\dag}$Department of Electrical and Computer Engineering,
The University of Texas at Austin,
EER~7.824, Austin, TX~~78712}
\email{ari@utexas.edu}

\author[Anup Biswas]{Anup Biswas$^\ddag$}
\address{$^\ddag$Department of Mathematics,
Indian Institute of Science Education and Research,
Dr.\ Homi Bhabha Road, Pune 411008, India}
\email{anup@iiserpune.ac.in}

\author[Somnath Pradhan]{Somnath Pradhan$^\ddag$}
\email{somnath@iiserpune.ac.in}

%%%%%%%%%%%%%%%%%%%%%%%%%%%%%%%%%%%%%%%%%%%%%%%%%%%%%%%%%%%%%%%%%%%%%%%%%%%%%%%%
\begin{abstract}
In this article we consider the ergodic risk-sensitive control problem
for a large class of multidimensional controlled diffusions on the
whole space.
We study the minimization and maximization problems under
either a blanket stability hypothesis,
or a near-monotone assumption on the
running cost.
We establish the convergence of
the policy improvement algorithm for these models.
We also present a more general
result concerning the region of attraction
of the equilibrium of the algorithm.
\end{abstract}
%%%%%%%%%%%%%%%%%%%%%%%%%%%%%%%%%%%%%%%%%%%%%%%%%%%%%%%%%%%%%%%%%%%%%%%%%%%%%%%%
\keywords{Principal eigenvalue, semilinear differential equations,
stochastic representation, policy improvement}

\subjclass[2000]{Primary 35P30, 93E20, 60J60}

%%%%%%%%%%%%%%%%%%%%%%%%%%%%%%%%%%%%%%%%%%%%%%%%%%%%%%%%%%%%%%%%%%%%%%%%%%%%%%%%
\maketitle

%%%%%%%%%%%%%%%%%%%%%%%%%%%%%%%%%%%%%%%%%%%%%%%%%%%%%%%%%%%%%%%%%%%%%%%%%%%%%%%%

\section{Introduction}

Given controlled dynamics governed by the It\^o equation
$$\D X_t \,=\, b(X_t,U_t)\,\D{t} + \upsigma(X_t)\,\D W_t$$
for some suitable functions $b$, and $\upsigma$, where $U$ is an admissible control
taking values in a compact metric space $\Act$, and $W$ is a standard
Wiener process, we consider the problem of minimizing or maximizing
the functional
\begin{equation*}
\sE_x(c, U) \,\df\,
\limsup_{T\to\infty} \, \frac{1}{T}\,
\log \Exp^U_x\Bigl[\E^{\int_0^T c(X_s, U_s)\,\D{s}}\Bigr]
\end{equation*}
over all admissible controls $U$, where $c$ is a
suitable \emph{running cost} function.
This is of course known as the (ergodic) risk-sensitive control problem.
The presence of the exponential in the functional reduces the possibility of rare but
devastating large excursions of the state process.
Risk-sensitive control has attracted a lot of interest because of its applications in
large deviations \cite{Kaise-06}, mathematical finance \cite{BP99,FS00},
missile guidance \cite{Speyer}, and other fields.
For a book-length exposition of this topic see \cite{Whittle}.
Though this topic has been around for quite sometime, some of major
issues such as uniqueness of value functions, verification results,
variational representations etc., for the (ergodic) risk-sensitive control
problem for a controlled diffusion were 
resolved fairly recently \cite{AB18,ABS19,AB18a,ABBK-19}.
The goal of this article is to establish the convergence
of the policy improvement algorithm (PIA) for this problem.
We are interested in both the minimization and maximization problems.
Due to the nonlinear nature of the criterion (w.r.t. $c$), these two problems
behave quite differently.
This article complements the results of \cite{BorMey02},
where a policy improvement algorithm is considered for skip-free countable state 
controlled Markov chains with a finite action set, under a near-monotone
cost assumption.

To describe our methodology, let us consider the minimization problem 
which seeks to determine 
$$\sEmin\,\df\, \inf_{x\in\Rd}\, \inf_{U\in\Uadm}\, \sE_x(c, U)\,,$$
where $\Uadm$ denotes the set of all admissible control.
As well known, under suitable hypotheses, $\sEmin$ is the principal eigenvalue of
a certain semilinear PDE.
More precisely, there exists a positive $V\in\cC^2(\Rd)$ satisfying
$$\trace\bigl(a(x)\grad^2 V(x)\bigr) + \min_{\zeta\in\Act}\,
\bigl[ b(x,\zeta)\cdot\grad V(x) + c(x,\zeta) V(x)\bigr]\,=\,\sEmin V(x)\,,$$
where $a\df\frac{1}{2}\,\upsigma\upsigma\transp$.
Existence, uniqueness and verification of optimality of stationary Markov controls
are studied in \cite{ABS19, AB18}.
For an early account of this problem we refer to
\cite{FM95,Menaldi-05,Nagai-96,Biswas-11,Biswas-11a}.
The PIA
can be described as follows: (1) we start with some stationary Markov control $v_0$,
(2) we determine the principal eigenpair $(\lambda_k,V_k)$, $k\geq 0$, solving
the eigenvalue problem in $\Rd$ with drift $b(x,v_k)$ and running
cost $c(x,v_k)$, and (3) we improve the control by selecting
$v_{k+1}$ from 
$$\Argmin_{\zeta\in\Act}\,\bigl[ b(x,\zeta)\cdot\grad V_k(x)
+ c(x,\zeta) V_k(x)\bigr]\,.$$
We say that the PIA converges if $(\lambda_k, V_k)$ converges to
$(\sEmin, V)$, in a suitable sense, as $k\to\infty$.
Similar problems have been studied for the
ergodic control problem with an additive criterion;
see for instance \cite{Krylov,Ari-12} and references therein.
Recently, \cite{AMV-17,AMV-18} study the PIA for ergodic control problems
in dimension $1$, where the diffusion matrix is also allowed to depend on the
control parameter.
The analysis of the PIA for the risk-sensitive problem is very different from the
the ergodic control problems with additive criteria.
Our analysis heavily uses the monotonicity property of the principal eigenvalue,
and the recurrence properties of the \emph{ground state process}
\cite{ABS19}.
The main results are \cref{T3.2} in \cref{S-minrisk}
which studies the minimization problem
under uniform ergodicity hypotheses (see \cref{A2.1}), and \cref{T4.2}
in \cref{S-maxrisk} which deals with
maximization problem under the near-monotone hypothesis in \cref{A2.2}.
Also, \cref{T5.2} in \cref{S5} which
is devoted to a more abstract study of the convergence of the PIA,
and finally \cref{T5.4} for the minimization problem under
a near-monotone hypothesis on the running cost.

%%%%%%%%%%%%%%%%%%%%%%%%%%%%%%%%%%%%%%%%%%%%%%%%%%%%%%%%%%%%%%%%%%%%%%%%%%%%%%%%
\subsection{Notation}\label{Snot}
We denote by $\uptau(A)$ the \emph{first exit time} of the process
$\{X_{t}\}$ from the set $A\subset\RR^{d}$, defined by
\begin{equation*}
\uptau(A) \,\df\, \inf\,\{t>0\,\colon X_{t}\not\in A\}\,.
\end{equation*}
The open ball of radius $r$ centered at $x\in\Rd$
is denoted by $\sB_{r}(x)$, and $\sB_{r}$ without an argument denotes
the ball centered at $0$.
We let $\uptau_{r}\df \uptau(\sB_{r})$,
and $\uuptau_{r}\df \uptau(\sB^{c}_{r})$.

The complement and closure
of a set $A\subset\Rd$ are denoted
by $A^{c}$ and $\Bar{A}$, respectively, and  $\Ind_A$
denotes its indicator function.
Given $a,b\in\RR$, the minimum (maximum) is denoted by $a\wedge b$ 
($a\vee b$), respectively,
and $a^\pm \df (\pm a)\vee 0$.
The inner product of two vectors $x$ and $y$ in $\Rd$ is denoted
as $x\cdot y$, or $\langle x,y\rangle$, $\abs{\,\cdot\,}$ denotes
the Euclidean norm, $x\transp$ stands for
the transpose of $x$, and
$\trace S$ denotes the trace of a square matrix $S$.

The term \emph{domain} in $\RR^{d}$
refers to a nonempty, connected open subset of the Euclidean space $\RR^{d}$. 
For a domain $D\subset\RR^{d}$,
the space $\cC^{k}(D)$ ($\cC^{k}_{b}(D)$) $k\ge 0$,
refers to the class of all real-valued functions on $D$ whose partial
derivatives up to order $k$ exist and are continuous (and bounded),
$\cC_{\mathrm{c}}^k(D)$ denotes its subset
 consisting of functions that have compact support,
and $\cC_0^k(D)$ the closure of $\cC_{\mathrm{c}}^k(D)$.
The space $\Lp^{p}(D)$, $p\in[1,\infty)$, stands for the Banach space
of (equivalence classes of) measurable functions $f$ satisfying
$\int_{D} \abs{f(x)}^{p}\,\D{x}<\infty$, and $\Lp^{\infty}(D)$ is the
Banach space of functions that are essentially bounded in $D$.
The standard Sobolev space of functions on $D$ whose generalized
derivatives up to order $k$ are in $\Lp^{p}(D)$, equipped with its natural
norm, is denoted by $\Sob^{k,p}(D)$, $k\ge0$, $p\ge1$.
In general, if $\mathcal{X}$ is a space of real-valued functions on $Q$,
$\mathcal{X}_{\mathrm{loc}}$ consists of all functions $f$ such that
$f\varphi\in\mathcal{X}$ for every $\varphi\in\cC_{\mathrm{c}}(Q)$.
Likewise, we define $\Sobl^{k, p}(D)$.

The constants $\kappa_1, \kappa_2, \ldots$ are generic and their
values might differ from place to place.

%%%%%%%%%%%%%%%%%%%%%%%%%%%%%%%%%%%%%%%%%%%%%%%%%%%%%%%%%%%%%%%%%%%%%%%%%%%%%%%%
\section{Model and Assumptions}
The controlled diffusion process $\process{X}$ in $\Rd$ is
governed by the It\^{o} equation
\begin{equation}\label{E2.1}
\D X_t \,=\, b(X_t,U_t)\,\D{t} + \upsigma(X_t)\,\D W_t\,,
\quad X_0=x\in\Rd\,.
\end{equation}
Here, $W$ is a $d$-dimensional standard Wiener process defined on a complete
probability space $(\Omega, \sF, \Prob)$.
The control process $\process{U}$ takes values in a compact metric space $\Act$,
is progressively measurable with respect to $\sF_t$, and is \emph{non-anticipative}:
for $s<t$, $W_t-W_s$ is independent of
\begin{equation*}
\sF_s \,\df\, \text{the\ completion\ of\ }
\sigma\{X_0, U_r, W_r  \,\colon\, r\le s\}
\text{\ relative\ to\ } (\sF, \Prob)\,.
\end{equation*}
The process $U$ is called an \emph{admissible} control,
and the set of all admissible control is denoted by $\Uadm$. 

Let $a= \frac{1}{2}\upsigma \upsigma\transp$.
We impose the following assumptions to guarantee existence of solution of \cref{E2.1}.

%%%%%%%%%%%%%%%%%%%%%%%%%%%%%%%%%%%%%%%%%%%%%%%%%%%%%%%%%%%%%%%%%%%%%%%%%%%%%%%%
\begin{enumerate}[(A2)]
\item[\hypertarget{A1}{(A1)}]
\emph{Local Lipschitz continuity:\/}
for some constant $C_{R}>0$ depending on $R>0$, the function
$b\colon\Rd\times\Act\to\Rd$, satisfies
\begin{equation*}
\abs{b(x,\zeta)-b(y, \zeta)}^2 
\,\le\, C_{R}\,\abs{x-y}^2
\end{equation*}
for all $x,y\in \sB_R$, and $\zeta\in\Act$. Furthermore, $\upsigma$ is assumed
to be bounded and globally Lipschitz.

\smallskip
\item[\hypertarget{A2}{(A2)}]
\emph{Growth condition:\/}
For some constant $C_0>0$, we have
\begin{equation}\label{EA2}
\sup_{\zeta\in\Act}\, \abs{b(x,\zeta)}  \,\le\,C_0
\bigl(1 + \abs{x}\bigr) \qquad \forall\, x\in\RR^{d}\,.
\end{equation}

\smallskip
\item[\hypertarget{A3}{(A3)}]
\emph{Nondegeneracy:\/}
For some positive constant $C$, it holds that
\begin{equation*}
\sum_{i,j=1}^{d} a^{ij}(x)\eta_{i}\eta_{j}
\,\ge\,C \abs{\eta}^{2} \qquad\forall\, x\in \Rd\,,\ \forall\,
\eta=(\eta_{1},\dotsc,\eta_{d})\transp\in\Rd\,.
\end{equation*}
\end{enumerate}

It is well known that under hypotheses
\hyperlink{A1}{(A1)}--\hyperlink{A2}{(A2)},
the stochastic differential equation in
\cref{E2.1} has a unique strong solution for
every admissible control (see for example, \cite[Theorem~ 2.2.4]{book}).
By a Markov control, we mean an admissible control of the form $v(t,X_t)$
for some Borel measurable function $v\colon \RR_+\times\Rd\to\Act$.
If $v$ is independent of $t$, we call it a stationary Markov control, and the set
of all stationary Markov controls is denoted by $\Usm$.
The hypotheses in \hyperlink{A1}{(A1)}--\hyperlink{A3}{(A3)}
imply the existence of unique
strong solutions under Markov controls \cite[Theorem~2.8]{Gyongy-96}.
We also mention here that for bounded measurable coefficients
Veretennikov  \cite{Veretennikov80} established
the existence of a strong solution by proving the pathwise uniqueness
(see also Theorem~5 in \cite{Veretennikov-82}).

%%%%%%%%%%%%%%%%%%%%%%%%%%%%%%%%%%%%%%%%%%%%%%%%%%%%%%%%%%%%%%%%%%%%%%%%%%%%%%%%
\begin{definition}
Let $c\colon\Rd\times\Act\to\RR$ be a continuous function
which represents the \emph{running cost}.
We assume that $c$ is locally Lipschitz continuous in 
its first argument uniformly with respect to the second.
Given a control $U\in\Uadm$ the \emph{risk-sensitive criterion} is given by
\begin{equation*}
\sE_x(c, U) \,\df\,
\limsup_{T\to\infty} \, \frac{1}{T}\,
\log \Exp^U_x\Bigl[\E^{\int_0^T c(X_s, U_s)\,\D{s}}\Bigr]\,,
\end{equation*}
and the optimal values are defined as
\begin{equation*}
\sEmin \,\df\, \inf_{x\in\Rd}\, \inf_{U\in\Uadm}\, \sE_x(c,U)\,,
\end{equation*}
and
\begin{equation*}
\sEmax \,\df\, \sup_{x\in\Rd}\, \sup_{U\in\Uadm}\, \sE_x(c,U)\,.
\end{equation*}
\end{definition}

These optimal values are extensively studied in 
\cite{ABBK-19,ABS19,AB18,AB18a,Biswas-11,Biswas-11a,FM95,Menaldi-05}.
In this article, we impose the following assumption on
the growth of $c$.
\begin{enumerate}[(A4)]
\item[\hypertarget{A4}{(A4)}]
\emph{Growth of the running cost:\/} For some constant $C$, we have 
\begin{equation*}
\sup_{\zeta\in\Act}\, \abs{c(x,\zeta)}  \,\le\,C
\bigl(1 + \abs{x}^2\bigr) \qquad \forall\, x\in\RR^{d}\,.
\end{equation*}
\end{enumerate}

Hypotheses \hyperlink{A1}{(A1)}--\hyperlink{A4}{(A4)} are assumed throughout
the paper, unless explicitly indicated otherwise.

%%%%%%%%%%%%%%%%%%%%%%%%%%%%%%%%%%%%%%%%%%%%%%%%%%%%%%%%%%%%%%%%%%%%%%%%%%%%%%%%
\begin{definition}
We define the operators $\cA$ and $\cA^c$ mapping
$\cC^2(\Rd)$ to $\cC(\Rd\times\Act)$ by
\begin{equation*}
\begin{aligned}
\cA u(x,\zeta) &\,\df\, \trace\bigl(a(x)\grad^2 u(x)\bigr) + b(x,\zeta)\cdot\grad u(x)\,,
\\[3pt]
\cA^c u(x,\zeta) &\,\df\, \cA u(x,\zeta) + c(x,\zeta) u(x)\,,
\end{aligned}
\end{equation*}
and for $v\in\Usm$, we often use the simplifying notation
\begin{equation}\label{Esimpl}
b_v(x) \,\df\, b\bigl(x,v(x)\bigr)\,,\quad\text{and\ \ }
c_v(x) \,\df\, c\bigl(x,v(x)\bigr)\,.
\end{equation}
\end{definition}

In the study of $\sEmin$ we enforce the following Foster--Lyapunov
condition on the dynamics.

%%%%%%%%%%%%%%%%%%%%%%%%%%%%%%%%%%%%%%%%%%%%%%%%%%%%%%%%%%%%%%%%%%%%%%%%%%%%%%%%
\begin{assumption}\label{A2.1}
In (a) and (b) below, $\Lyap \in \cC^2(\Rd)$ is some function taking
values in $[1,\infty)$, $c$ is nonnegative, $\widehat{C}$ is a positive constant,
and $\cK\subset\Rd$ is a compact set.
\begin{enumerate}
\item[(a)]
If $c$ is bounded, we assume without loss of generality
that $\inf_{\Rd\times\Act}\,c=0$, and that there exists some
constant $\gamma>\norm{c}_\infty$ satisfying
\begin{equation}\label{EA2.1A}
\cA \Lyap(x,\zeta)
\,\le\, \widehat{C} \Ind_{\cK}(x)-\gamma \Lyap(x)\qquad \forall\,(x,\zeta)
\in\Rd\times\Act\,.
\end{equation}

\item[(b)] If $c$ is not bounded, we assume
that there exists an inf-compact function $\ell$
(i.e., the sublevel sets $\{\ell \le \kappa\}$ are compact, or empty,
in $\Rd$ for each $\kappa\in\RR$)
such that $x\mapsto\ell(x)-\max_{\zeta\in\Act} c(x, \zeta)$ is also inf-compact, and
\begin{equation}\label{EA2.1B}
\cA \Lyap(x,\zeta)
\,\le\, \widehat{C} \Ind_{\cK}(x)-\ell(x) \Lyap(x)
\qquad \forall\,(x,\zeta)\in\Rd\times\Act\,.
\end{equation}
\end{enumerate}
\end{assumption}

As well known (see \cite{ABS19}), if $a$ and $b$ are bounded,
it might not be possible to find an unbounded function $\ell$
satisfying \cref{EA2.1B}.
This is the reason for \cref{EA2.1A}. Also, due to \cref{EA2}, it is known from
\cite[pg.~65]{book} that $x\mapsto\Exp^U_x[\uptau(\cK^c)]$ is inf-compact
for any 
admissible control $U$, and therefore, the Lyapunov function $\Lyap$ in
\cref{EA2.1A,EA2.1B} are also inf-compact.

Before we proceed further, let us
exhibit a class of dynamics satisfying \cref{A2.1}.

%%%%%%%%%%%%%%%%%%%%%%%%%%%%%%%%%%%%%%%%%%%%%%%%%%%%%%%%%%%%%%%%%%%%%%%%%%%%%%%%
\begin{example}
Suppose that
$\sup_{\zeta\in\Act} b(x,\zeta)\cdot x\le - \kappa \abs{x}^\alpha$
outside a compact set
for some $\alpha\in [1, 2]$, and $a$ is bounded. 
Define $\Lyap (x)\,\df\,\exp(\delta \sqrt{\abs{x}^2+1})$.
Then an easy calculation shows that
\begin{equation*}
\begin{aligned}
\cA \Lyap(x) &\,\le\, \kappa_1 \biggl(\Ind_{\cK_1}(x)
+ \delta\frac{1}{\sqrt{\abs{x}^2+1}}
+ \delta^2 \frac{\abs{x}^2}{\abs{x}^2+1}\biggr)\Lyap(x)
- \delta\frac{\abs{x}^\alpha}{\sqrt{\abs{x}^2+1}}\Lyap(x)
\end{aligned}
\end{equation*}
for some constant $\kappa_1$, and a compact set $\cK_1$.
Thus, if $\alpha>1$,
and we choose $\ell\sim\abs{x}^{\alpha-1}$, \cref{EA2.1B} is satisfied.
For $\alpha=1$,  choosing $\delta$ suitably small we obtain \cref{EA2.1A}.
\end{example}

For the maximization problem, we use the near-monotone hypothesis in
\cref{A2.2} below,
which is somewhat weaker than \cite[Assumption~3.1\,(iii)]{ABBK-19}.

%%%%%%%%%%%%%%%%%%%%%%%%%%%%%%%%%%%%%%%%%%%%%%%%%%%%%%%%%%%%%%%%%%%%%%%%%%%%%%%%
\begin{assumption}\label{A2.2}
It holds that 
\begin{equation}\label{EA2.2}
\lim_{r\to\infty}\,\max_{(x,\zeta)\in\sB_r^c\times\Act} c(x,\zeta)
\,<\,\sEmax \,<\, \infty\,.
\end{equation}
\end{assumption}

We let $\varrho_*$ denote the principal eigenvalue
of the maximal operator,
defined as follows
$$\varrho_*\,\df\,\inf\,\bigl\{\lambda\,\colon \exists\,\psi\in \Sobl^{2,p}(\Rd),\,
p>d,\, \psi>0, \text{\ satisfying\ }\max_{\zeta\in\Act}\cA(x,\zeta)\psi\le \lambda\psi
 \text{\ in\ } \Rd\bigr\}\,.$$
\cref{A2.2} implies that $\varrho_*$ is finite.
To see this, note first that
$\varrho_*>-\infty$, since the Dirichlet eigenvalue on the
unit ball is finite.
Also, there exists some $r_0\in(0,\infty)$ such that
\begin{equation*}
\max_{(x,\zeta)\in\sB_{r_0}^c\times\Act} c(x,\zeta) \,<\, \sEmax
\end{equation*}
by \cref{EA2.2}.
Hence, we have
\begin{equation*}
\begin{aligned}
\varrho_* &\,\le\, \max_{(x,\zeta)\in \Rd\times\Act}\, c(x,\zeta)\\
&\,<\, 
\max_{(x,\zeta)\in \sB_{r_0}^c\times\Act}\,  c(x,\zeta)
+ \max_{(x,\zeta)\in \sB_{r_0}\times\Act}\,  c(x,\zeta)
\\
&\,<\,\sEmax + \max_{(x,\zeta)\in \sB_{r_0}\times\Act}\,  c(x,\zeta)
\,<\, \infty\,.
\end{aligned}
\end{equation*}

%%%%%%%%%%%%%%%%%%%%%%%%%%%%%%%%%%%%%%%%%%%%%%%%%%%%%%%%%%%%%%%%%%%%%%%%%%%%%%%%
\begin{remark}
One can also study a maximization problem under \cref{A2.1}; the results and
proofs
are similar to the minimization problem. Therefore, we do not discuss the maximization
problem under \cref{A2.1}.
\end{remark}

%%%%%%%%%%%%%%%%%%%%%%%%%%%%%%%%%%%%%%%%%%%%%%%%%%%%%%%%%%%%%%%%%%%%%%%%%%%%%%%%
\subsection{Principal eigenvalues of linear operators}

In this section we recall some facts about principal eigenvalues
which are needed later.
Let $b, f\colon\Rd\to\RR$ be locally bounded Borel measurable functions.
We also 
assume that $\langle x, b(x)\rangle^+\le C(1+\abs{x}^2)$ for $x\in\Rd$.
Consider the linear operator
$$\sL^f u(x) \,\df\, \trace\bigl(a(x)\grad^2 u(x)\bigr)
+ b(x)\cdot\grad u(x) + f(x) u(x)\,.$$
The principal eigenvalue $\lamstr(\sL^f)$ of $\sL^f$ is defined
as follows \cite{Berestycki-15}:
\begin{equation}\label{E-lamstr}
\lamstr(\sL^f) \,\df\, \inf\,\bigl\{\lambda\in\RR\,\colon \exists\,\text{positive\ }
\psi\in\Sobl^{2,p}(\Rd)\,,\; p>d\,, \text{\ satisfying\ }
 \sL^f\psi\le\lambda\psi \text{\ in\ } \Rd\bigl\}\,.
\end{equation}
The reader should have noticed that the convention in this definition is
the opposite of what is normally used in the pde literature.
See for example the definition in \cite[eq.~(1.10)]{Berestycki-94}.
The same convention is used in \cite{AB18,AB18a,ABS19} which contain results that
we cite in this paper.

When we want to emphasize the dependence of $\lamstr(\sL^f)$ on $f$,
we write this as $\lamstr(f)$.
It is also known from \cite[Theorem~1.4]{Berestycki-15}, that for any
$\lambda\in[\lamstr(f), \infty)$, there exists a positive $\Psi\in \Sobl^{2,p}(\Rd)$,
$p>d$, satisfying $\sL^f\Psi=\lambda\Psi$ in $\Rd$.
We denote the principal eigenpair by $(\lamstr,\Psi^*)$.
We now recall the following definition of 
strict monotonicity from \cite{ABS19}

%%%%%%%%%%%%%%%%%%%%%%%%%%%%%%%%%%%%%%%%%%%%%%%%%%%%%%%%%%%%%%%%%%%%%%%%%%%%%%%%
\begin{definition}\label{D2.3}
We say that $\lamstr$ is \emph{strictly monotone at} $f$, if for some non-trivial
nonnegative function $h$, vanishing at infinity, we have
$\lamstr(f-h)<\lamstr(f)$.
We also say that $\lamstr$ is \emph{strictly monotone at $f$ on the right},
if for any non-trivial
nonnegative function $h$, we have
$\lamstr(f)<\lamstr(f+h)$.
\end{definition}

Note that strictly monotonicity, implies strictly monotonicity on the right,
since the eigenvalue is a convex function of $f$.

Given an eigenpair $(\lambda,\Psi)$ we define
the \textit{twisted process} $\process{Y}$ as a solution to the SDE
\begin{equation}\label{E-twist}
\D{Y_t} \,=\, b(Y_t)\,\D{t} + 2a(Y_t) \grad \psi(Y_t)\,\D{t} + \upsigma(Y_t)\,\D{W_t}\,,
\end{equation}
where $\psi=\log\Psi$.
The twisted process corresponding to a principal eigenpair $(\lamstr,\Psi^*)$
is called a \textit{ground state process}.
Then we have the following result from \cite[Theorem~2.1]{ABS19}.

%%%%%%%%%%%%%%%%%%%%%%%%%%%%%%%%%%%%%%%%%%%%%%%%%%%%%%%%%%%%%%%%%%%%%%%%%%%%%%%%
\begin{theorem}\label{T2.1}
Suppose that $\inf_{\Rd} f>-\infty$ and $\lamstr(\sL^f)$ is finite.
Then, the following hold:
\begin{enumerate}
\item[\ttup{i}] For $\lambda>\lamstr$, the twisted process \cref{E-twist}
corresponding to any eigenpair $(\lambda,\Psi)$ is transient.
\item[\ttup{ii}]
A ground state process is exponentially ergodic if and only if
$\lamstr$ is strictly monotone at $f$.
\end{enumerate}
\end{theorem}

%%%%%%%%%%%%%%%%%%%%%%%%%%%%%%%%%%%%%%%%%%%%%%%%%%%%%%%%%%%%%%%%%%%%%%%%%%%%%%%%
\begin{remark}\label{R2.2}
The results of \cref{T2.1} also hold if $\sup_{\Rd} f<\infty$.
\end{remark}

%%%%%%%%%%%%%%%%%%%%%%%%%%%%%%%%%%%%%%%%%%%%%%%%%%%%%%%%%%%%%%%%%%%%%%%%%%%%%%%%
\section{Policy improvement for the minimization problem}\label{S-minrisk}
In this section we prove the convergence
of the
policy improvement algorithm (PIA) for $\sEmin$ under \cref{A2.1}.
We begin with the following result which is an extension of 
\cite[Theorem~4.1]{ABS19}.

%%%%%%%%%%%%%%%%%%%%%%%%%%%%%%%%%%%%%%%%%%%%%%%%%%%%%%%%%%%%%%%%%%%%%%%%%%%%%%%%
\begin{theorem}\label{T3.1}
Grant \cref{A2.1}. There exists a positive solution $V\in\cC^2(\Rd)$ satisfying
\begin{equation}\label{ET3.1A}
\min_{\zeta\in\Act}\, \cA^c V(x,\zeta) \,=\,
 \sEmin V(x)\quad \text{in\ } \Rd, \quad \text{and}\quad V(0)=1\,.
\end{equation}
In addition, with $\bUsm\subset\Usm$ denoting
the class of Markov controls $v$ which satisfy
\begin{equation*}
b\bigl(x,v(x)\bigr)\cdot\grad{V}(x) + c\bigl(x,v(x)\bigr){V}(x)\,=\,
\min_{\zeta\in\Act}\,\bigl[ b(x,\zeta)\cdot\grad{V}(x) + c(x,\zeta){V}(x)\bigr]
\quad \text{a.e.\ in}\ \Rd\,,
\end{equation*}
the following hold:
\begin{itemize}
\item[(a)] $\bUsm\subset \Usms$ and $\lamstr(c_v)=\sEmin$ for all $v\in\bUsm$.
Here, $\Usms$ denotes the set of all optimal stationary Markov controls.
\item[(b)] $\Usms\subset \bUsm$.
\item[(c)] \Cref{ET3.1A} has a unique solution in $\cC^2(\Rd)$.
\end{itemize}
\end{theorem}

Next, we describe the PIA.

%%%%%%%%%%%%%%%%%%%%%%%%%%%%%%%%%%%%%%%%%%%%%%%%%%%%%%%%%%%%%%%%%%%%%%%%%%%%%%%%
\begin{algorithm}\label{Alg3.1}
Policy iteration.
\begin{itemize}
\item[1.] Initialization. Set $k=0$ and select any $v_0\in\Usm$.
\smallskip
\item[2.] Value determination. Let $V_k\in \Sobl^{2,p}(\Rd)$, $p>d$,
be the unique principal eigenfunction satisfying $V_k(0)=1$ and
\begin{equation}\label{E3.2}
\trace\bigl(a(x)\grad^2 V_k(x)\bigr) + b(x,v_k(x))\cdot\grad V_k(x) + c(x, v_k(x)) V_k(x)
\,=\, \lamstr(c_{v_k}) V_k(x)\,,\quad x\in\Rd\,.
\end{equation}
Existence of a unique principal eigenfunction in \cref{E3.2} follows from
\cite[Section~3]{ABS19}. We let $\lambda_k\df\lamstr(c_{v_k})$.
\smallskip
\item[3.] Policy improvement. Choose any $v_{k+1}\in\Usm$ satisfying
$$v_{k+1}(x)\in \Argmin_{\zeta\in\Act}\,
\left[b(x,\zeta)\cdot\grad V_k(x) + c(x,\zeta) V_k(x)\right], \quad x\in\Rd\,.$$
\end{itemize}
\end{algorithm}

The main result of this section is the following.

%%%%%%%%%%%%%%%%%%%%%%%%%%%%%%%%%%%%%%%%%%%%%%%%%%%%%%%%%%%%%%%%%%%%%%%%%%%%%%%%
\begin{theorem}\label{T3.2}
Under \cref{A2.1}, the following hold:
\begin{itemize}
\item[(i)] For all $k\in\NN$, we have $\lambda_{k+1}\le \lambda_{k}$,
and $\lim_{k\to\infty}\lambda_k=\sEmin$.
\item[(ii)] The sequence $\{V_k\}$ converges weakly in $\Sobl^{2,p}(\Rd)$, $p>d$,
to the unique solution $V$ of \cref{ET3.1A}.
\end{itemize}
\end{theorem}

In the sequel, we use the notation
\begin{equation}\label{E-notaA}
\begin{gathered}
b_k(x) \,\df\, b\bigl(x,v_k(x)\bigr)\,,\qquad
c_k(x) \,\df\, c\bigl(x,v_k(x)\bigr)\,,\\
\sL_k f(x)\,\df\, \trace\bigl(a(x)\grad^2 f(x)\bigr)
+ b_k(x)\cdot\grad f(x) + c_k(x) f(x)\,.
\end{gathered}
\end{equation}
Let us start with following lemma.

%%%%%%%%%%%%%%%%%%%%%%%%%%%%%%%%%%%%%%%%%%%%%%%%%%%%%%%%%%%%%%%%%%%%%%%%%%%%%%%%
\begin{lemma}\label{L3.1}
We have $\lambda_{k+1}\le \lambda_k$ for all $k\ge 1$.
\end{lemma}
\begin{proof}
From the policy improvement algorithm, we have
\begin{equation}\label{EL3.1A}
\begin{aligned}
\cA V_{k}\bigl(x,v_{k+1}(x)\bigr) + c\bigl(x,v_{k+1}(x)\bigr) V_{k}(x)
& \,=\, \min_{\zeta\in\Act}\, \cA^c V_{k}(x,\zeta)\\
& \,\le\, \cA V_{k}\bigl(x,v_{k}(x)\bigr) + c\bigl(x,v_{k}(x)\bigr) V_{k}(x)  \\ 
& \,= \, \lambda_{k} V_{k}(x)\quad \text{a.e.\ in\ } \Rd\,.
\end{aligned}
\end{equation}
It then follows from the definition of the principal eigenvalue that
$\lambda_{k+1}\le \lambda_{k}$.
\end{proof}

Next we recall the following global estimate
from \cite[Lemma~4.1]{ABBK-19}.

%%%%%%%%%%%%%%%%%%%%%%%%%%%%%%%%%%%%%%%%%%%%%%%%%%%%%%%%%%%%%%%%%%%%%%%%%%%%%%%%
\begin{lemma}\label{L3.2}
Suppose that $a$ is bounded and has a uniform modulus of continuity in $\Rd$,
and
$b\colon\Rd\to\Rd$, $c\colon\Rd\to\RR$ are locally bounded. Then
there exists a constant $\widetilde{C}$, dependent on $a$, such that for any
strong positive solution $\phi\in\Sobl^{2,p}(\Rd)$, $p>d$, of $\cA^c\phi=0$ we have
\begin{equation*}
\frac{\abs{\grad \phi(x)}}{\phi(x)}\,\le\, \widetilde{C}\left(1+\sup_{y\in\sB_1(x)}
\left(|b(y)| + \sqrt{|c(y)|}\right)\right), \quad x\in\Rd\,.
\end{equation*}
\end{lemma}

Continuing, we define
\begin{equation}\label{E-psi}
\begin{aligned}
\psi_{k}(x) &\,\df\, -\frac{1}{V_{k-1}(x)}\,\cA V_{k-1}\bigl(x,v_{k}(x)\bigr)
- c\bigl(x,v_{k}(x)\bigr)  + \lambda_{k-1}\\
&\,=\,
-\frac{1}{V_{k-1}(x)}
\min_{\zeta\in\Act}\, \cA^c V_{k-1}(x,\zeta) + \lambda_{k-1}\\
&\,=\,
\frac{1}{V_{k-1}(x)}\Bigl(
b\bigl(x,v_{k-1}(x)\bigr)\grad V_{k-1}(x)+c\bigl(x,v_{k-1}(x)\bigr)V_{k-1}(x)\\
&\mspace{150mu}
- \min_{\zeta\in\Act}\,\bigl[b(x,\zeta)\cdot\grad V_{k-1}(x)
+ c(x,\zeta) V_{k-1}(x)\bigr]\Bigr)
\,,\quad k\in\NN\,.
\end{aligned}
\end{equation}
Note that $\psi_{k}$ is a nonnegative function.
We write \cref{EL3.1A} in the form
\begin{equation}\label{E3.5}
\min_{\zeta\in\Act}\, \cA^c V_{k} (x,\zeta) - \lambda_{k}V_{k}(x) 
\,=\,  - \psi_{k+1}(x)V_{k}(x) \,\le\, 0 \quad \text{a.e.\ in\ } \Rd\,.
\end{equation}
Applying \cref{L3.2} to the equation 
\begin{equation}\label{E3.6}
\sL_k V_{k}(x)\,=\, \lambda_{k} V_{k}(x)\quad \text{a.e.\ in\ } \Rd,
\end{equation} 
we get that 
\begin{equation}\label{E3.7}
\frac{\abs{\grad V_k(x)}}{V_k(x)} \,\le\,
 \widetilde{C}\biggl(1+\sup_{y\in\sB_1(x)}\left(|b_k(y)| + 
\sqrt{|c_k(y)-\lambda_k|}\right)\biggr)
\,\le\, \widetilde{C}_1 (1+\abs{x})\,,\quad x\in\Rd,
\end{equation}
for some constant $\widetilde{C}_1$,
where in the last inequality we use \hyperlink{A2}{(A2)} and \hyperlink{A4}{(A4)}.
Again, using \cref{E3.6} together with \cref{E3.7}, we find that
\begin{align*}
\frac{1}{V_k(x)}
\babs{\trace\bigl(a(x)\grad^2 V_k(x)\bigr)}\,\le\,
|b_k(x)|\frac{|\grad V_k(x)|}{V_k(x)} + |c_k(x)-\lambda_k|
\,\le\, \widetilde{C}_2(1+\abs{x}^2)\,.
\end{align*}
Thus, from \cref{E3.5}, we obtain
\begin{equation}\label{E3.8}
\psi_{k+1}(x) \,\le\, \widetilde{C}_3(1 + |x|^2),\quad  \text{for all\ } x\in\Rd,
\end{equation}
for some constant $\widetilde{C}_3$.

Let $\process{Y^{k}}$ be the twisted process \cref{E-twist} corresponding to the principal
eigenpair $(\lambda_{k},V_{k})$, that is, $\process{Y^{k}}$ is the unique solution to the
following SDE
\begin{equation}\label{E3.9}
\D{Y_t ^{k}} \,=\, b_k(Y_t ^{k})\,\D{t} + 2a(Y_t ^{k}) \grad \log V_{k}(Y_t ^{k})\,\D{t}
+ \upsigma(Y_t ^{k})\,\D{W_t}\,,
\end{equation}
and let
\begin{equation}\label{E3.10}
\widetilde\sL_{k} f \,\df\, (\sL_k-c_k(x)) f(x) + 2a(x)\grad \log V_k(x)\cdot\grad f(x)
\end{equation} 
denote the extended generator of \cref{E3.9}.
It follows by \cref{E3.7} that the drift of \cref{E3.9} has at most linear growth,
and, therefore, $\process{Y^{k}}$ is non-explosive. 
An easy calculation shows that,
for any positive $\phi\in \Sobl^{2,p}(\Rd)$, we have the following identity
\begin{equation}\label{E3.11}
\widetilde\sL_k \left(\frac{\phi}{V_{k}}\right) \,=\,
\left(\frac{\cA \phi}{\phi} - \frac{\cA V_{k}}{V_{k}}\right)\,
\left(\frac{\phi}{V_{k}}\right)\quad\text{a.e.\ in\ }\Rd\,.
\end{equation}

We let $\widetilde{P}^k_t(x,\D{y})$ denote the transition probability of
$\process{Y^k}$.

%%%%%%%%%%%%%%%%%%%%%%%%%%%%%%%%%%%%%%%%%%%%%%%%%%%%%%%%%%%%%%%%%%%%%%%%%%%%%%%%
\begin{lemma}\label{L3.3}
Grant \cref{A2.1}.
Then the family of invariant
probability measures $\{\Tilde\upmu_k\}_{k\in\NN}$ corresponding
to the operators $\{\widetilde\sL_k\}_{k\in\NN}$ is tight.
In addition,
there exist positive constants
$\gamma_\circ$ and $C_{\gamma_\circ}$ such that
\begin{equation}\label{EL3.3A}
\bnorm{\widetilde{P}^k_t(x,\cdot\,)-\Tilde\upmu_k(\cdot)\,}_{\tv}\,\le\,
C_{\gamma_\circ} \frac{\Lyap(x)}{V_k(x)}\, \E^{-\gamma_\circ t}\qquad
\forall\, (t,x)\in\RR_+\times\Rd\,,\ 
\forall\, k\in\NN\,.
\end{equation}
\end{lemma}

\begin{proof}
Consider \cref{A2.1}\,(a), and let
$\epsilon_\circ\ge \gamma-\norm{c}_\infty$.
By \cref{E3.6,E3.11}, we have
\begin{equation}\label{PL3.3A}
\begin{aligned}
\widetilde\sL_k \biggl(\frac{\Lyap}{V_k}\biggr)(x) &\,\le\,
\Bigl(\widehat{C}\Lyap^{-1}(x) \Ind_{\sK}(x) - \gamma + c_k(x) 
 -\lambda_k\Bigr) \biggl(\frac{\Lyap}{V_k}\biggr)(x)\\
&\,\le\,
\Bigl(\widehat{C} \Ind_{\sK}(x) -  \epsilon_\circ
-\lambda_k\Bigr) \biggl(\frac{\Lyap}{V_k}\biggr)(x)\quad\forall\,k\in\NN\,,
\end{aligned}
\end{equation}
 since $\Lyap^{-1}\le 1$.
\Cref{PL3.3A} implies that
\begin{equation}\label{PL3.3B}
\int_\Rd \frac{\Lyap(x)}{V_k(x)}\,\Tilde\upmu_k(\D{x})\,\le\, M
\quad \forall\,k\in\NN\,,
\end{equation}
for some constant $M$.
Note that here $\frac{\Lyap}{V_k}$ is, in general, not in $\cC^2(\Rd)$,
so to derive \cref{PL3.3B}, we apply \cite[Lemma~3.1]{ACPZ18}.
It is well known that the linear growth bound in \cref{E3.7} implies that
$x\mapsto
\inf_{k\in\NN}\,\widetilde\Exp_x^k[\uptau(\sK^c)]$ is inf-compact.
Recall the definition $\uuptau_r=\uptau(\sB_r^c)$ in \cref{Snot}.
Let $\{X_t^k\}$ be the process satisfying \eqref{E2.1} with the Markov control $v_k$ and
$\Exp^k[\cdot]$ be the corresponding expectation operator.
It is straightforward to show, using the stochastic representation
\begin{equation*}
V_k(x) \,=\, \Exp^k_x\left[\E^{\int_0^{\uuptau_r}(c_k(X_t^k)-\lambda_k)\,\D{t}}\,
V_k\bigl(X_{\uuptau_r}^k\bigr)\right]
\end{equation*}
which holds for any $r>0$ by \cite[Theorem~3.2]{ABS19},
and the inequality
\begin{equation*}
\Lyap(x) \,\ge\, \Exp^k_x\left[\E^{\gamma\uptau(\sK^c)}
\, \Lyap \left(X_{\uptau(\sK^c)}^k\right)\right]\,,
\end{equation*}
which follows from \cref{EA2.1A},
that
$\inf_{(k,x)\in\NN\times\Rd}\,\frac{\Lyap(x)}{V_k(x)}>0$.
On the other hand, applying the It\^o--Krylov formula
to \cref{PL3.3A}, followed by Fatou's lemma and the Jensen inequality,
we obtain
\begin{equation*}
\frac{\Lyap(x)}{V_k(x)}\,\ge\,\epsilon_\circ\,
\widetilde\Exp_x^k\bigl[\uptau(\sK^c)\bigr]\,
\inf_{(k,x)\in\NN\times\Rd}\,\frac{\Lyap(x)}{V_k(x)}\,,
\end{equation*}
which implies that $x\mapsto\inf_{k\in\NN}\,\frac{\Lyap(x)}{V_k(x)}$ is
inf-compact.
Hence, tightness of $\{\Tilde\upmu_k\}_{k\in\NN}$ follows by \cref{PL3.3B}.
The proof under \cref{A2.1}\,(b) is identical.

It is well known that \cref{PL3.3A} implies that $\process{Y^k}$ is
exponentially ergodic.
It remains to show that the bound $\gamma_\circ$ on the rate convergence,
and the constant $C_{\gamma_\circ}$ can be chosen independently of $k\in\NN$.
For this, we employ the same exact argument as in the proof
of \cite[Theorem~2.1\,(b)]{AHP18}, which utilizes the estimates in
\cite[Theorem~2.3]{MeynT-94b}.
This completes the proof.
\end{proof}

%%%%%%%%%%%%%%%%%%%%%%%%%%%%%%%%%%%%%%%%%%%%%%%%%%%%%%%%%%%%%%%%%%%%%%%%%%%%%%%%
\begin{lemma}\label{L3.4}
Under \cref{A2.1},
$\psi_k$ converges to $0$ in $\Lpl^p(\Rd)$ for any $p\ge1$.
\end{lemma}

\begin{proof}
Applying the It\^o-Krylov formula to
$\sL_{k}V_{k-1} + (\psi_k- \lambda_{k-1})V_{k-1}=0$,
we obtain
\begin{equation}\label{PL3.4A}
\int_0^T \Exp_x^k\Bigl[\E^{\int_0^t (c_k(X_s^k) - \lambda_{k-1})\,\D{s}}\,
\psi_k(X_t^k)V_{k-1}(X_t^k)\Bigr]\,\D{t} \,\le\, V_{k-1}(x)
\quad\forall\,k\in\NN\,,\ \forall\,T>0\,.
\end{equation}
Let
\begin{equation}\label{PL3.4B}
h_k\,\df\,\psi_k\,\frac{V_{k-1}}{V_k}\,,\quad k\in\NN\,,
\end{equation}
and $\Tilde\lambda_k\df\lambda_{k-1}-\lambda_{k}$.
By \cite[Lemma~2.3]{ABS19}, we have
\begin{equation}\label{PL3.4C}
\begin{aligned}
\Exp_x^k\Bigl[\E^{\int_0^t (c_k(X_s^k) - \lambda_{k-1})\,\D{s}}\,
\psi_k(X_t^k)V_{k-1}(X_t^k)\Bigr]
&\,=\, \E^{-\Tilde\lambda_k t}\,V_k(x)\,\widetilde\Exp_x^{k}\bigl[h_k(Y_t^k)\bigr]\\
&\,=\, \E^{-\Tilde\lambda_k t}\,V_k(x)\,
\int_\Rd \widetilde{P}^k_t(x,\D{y}) h_k(y)\,.
\end{aligned}
\end{equation}
It is clear, by Sobolev imbedding, that the families of functions $\{V_k\}_{k\in\NN}$,
$\{V_k^{-1}\}_{k\in\NN}$, and 
$\{\grad V_k\}_{k\in\NN}$ are locally H\"older equicontinuous.
Hence, $\{\psi_k\}_{k\in\NN}$ and $\{h_k\}_{k\in\NN}$
defined in \cref{E-psi,PL3.4B}, respectively,
are locally bounded.

If $\Tilde\lambda_k=0$, then we have
$\Tilde\sL_k(\frac{V_{k-1}}{V_k})\leq 0$ by \eqref{E3.11}, and
then by \cref{L3.3} we obtain $V_{k-1}=V_k$ implying $\psi_k=0$.
So we suppose, without loss of generality, that $\Tilde\lambda_k>0$.
Let $D\subset\Rd$ be an arbitrary ball. By \cref{EL3.3A}, we have
\begin{equation}\label{PL3.4F}
\int_D \widetilde{P}^{k}_t(x,\D{y}) h_{k}(y)
\,\ge\, \int_{D} h_{k}(y)\,\Tilde\upmu_{k}(\D{y}) -
\bnorm{h_{k}}_{\infty,D}\,
C_{\gamma_\circ} \frac{\Lyap(x)}{V_{k}(x)}\, \E^{-\gamma_\circ t}\,,
\end{equation}
where $\norm{\,\cdot\,}_{\infty,D}$ denotes the infinity norm
of the restriction of a function to $D$.
Thus, by \cref{PL3.4A,PL3.4C,PL3.4F} we obtain
\begin{equation*}
\frac{1}{\Tilde\lambda_{k}} V_{k}(x)\int_{D} h_{k}(y)\,\Tilde\upmu_{k}(\D{y}) -
\frac{C_{\gamma_\circ}}{\gamma_\circ} \bnorm{h_{k}}_{\infty,D}\,\Lyap(x)
\,\le\, V_{k-1}(x)\quad\forall\,k\in\NN\,,
\end{equation*}
which implies that
\begin{equation}\label{PL3.4G}
\int_{D} h_{k}(y)\,\Tilde\upmu_{k}(\D{y})
\,\le\, \frac{\Tilde\lambda_{k}}{V_{k}(x)}\Bigl(
V_{k-1}(x)
+\frac{C_{\gamma_\circ}}{\gamma_\circ} \bnorm{h_{k}}_{\infty,D}\,\Lyap(x)\Bigr)
\quad\forall\,k\in\NN\,.
\end{equation}
Since $\{\Tilde\upmu_k\}_{k\in\NN}$ is a tight family, and $a$ is uniformly elliptic,
and the drifts of \cref{E3.10} are bounded in $D$ uniformly in $k\in\NN$,
the corresponding densities $\rho_k$ are bounded above and away from $0$
uniformly in $k\in\NN$ (see \cite[Lemma~3.2.4\,(b)]{book}).
Thus, the result follows by \cref{PL3.4G}, since, for any fixed $x$,
$\{V_k(x)\}_{k\in\NN}$,
$\{V_k^{-1}(x)\}_{k\in\NN}$, and $\bigl\{\bnorm{h_{k}}_{\infty,D}\bigr\}_{k\in\NN}$
are bounded, whereas $\Tilde\lambda_{k}\searrow0$ as $k\to\infty$.
This completes the proof.
\end{proof}

%%%%%%%%%%%%%%%%%%%%%%%%%%%%%%%%%%%%%%%%%%%%%%%%%%%%%%%%%%%%%%%%%%%%%%%%%%%%%%%%
Now we are ready to prove \cref{T3.2}.
\begin{proof}[Proof of \cref{T3.2}]
From \cref{E3.5} we have
\begin{equation}\label{ET3.2A}
\min_{\zeta\in\Act}\, \cA^c V_k(x,\zeta) + \psi_{k+1}(x) V_k(x) -\lambda_k V_k(x)
\,=\, \sL_{k+1} V_k(x) + \psi_{k+1}(x) V_k(x) -\lambda_k V_k(x) \,=\,0\,.
\end{equation}
Since, $V_k$ solves \cref{E3.2} and $V_k(0)=1$, and $\lambda_k$ converges by \cref{L3.1},
it follows from the Harnack inequality that $V_k$ is bounded.
Thus, applying the well-known a priori estimate
\cite[Lemma~5.3]{ChenWu},
it follows that $\{V_k\}_{k\ge 0}$ is 
locally bounded in $\Sob^{2,p}(\Rd)$, for $p>d$.
Therefore, employing the Cantor diagonal argument we can extract a subsequence
$\{V_{n_k}\}$ that
converges weakly to some
$\widetilde{V}\in \Sobl^{2,p}(\Rd)$, $p>d$.
Furthermore,  we have
$V_{n_k}\to \widetilde{V}$ in $\cC^{1}_{\rm loc}(\Rd)$ by Sobolev embedding,
and $\widetilde{V}\in\sorder(\Lyap)$ from the proof of \cref{L3.3}.
Let $\Lambda=\lim_{k\to\infty}\lambda_k\ge 0$.
Then passing to the limit in \cref{ET3.2A}
and using \cref{L3.4} we obtain that
\begin{equation}\label{ET3.2B}
\min_{\zeta\in\Act}\, \cA^c \widetilde{V}(x,\zeta)  -\Lambda \widetilde{V}(x)
\,=\, 0\quad \text{in\ } \Rd\,.
\end{equation}
It follows from standard elliptic regularity theory
(see \cite[Theorem~9.19]{GilTru}) that $\widetilde{V}\in\cC^2(\Rd)$.
More precisely, since
$\widetilde{V}\in\Sobl^{2, p}(\Rd)$, for any $p>d$,
it follows from Sobolev embedding theorem that
$\widetilde{V}\in \cC^{1, \alpha}_{\rm loc}(\Rd)$ for any $\alpha\in(0,1)$. Now
using the locally Lipschitz property of $b$ and $c$ we can write \eqref{ET3.2B} as 
$$\trace\bigl(a(x)\grad^2 \widetilde{V}(x)\bigr)=f(x)$$
for some locally $\alpha$-H\"{o}lder continuous function $f$.
Applying \cite[Theorem~9.19]{GilTru}, we
then obtain that $\widetilde{V}\in\cC^{2,\alpha}(\sB_R)$ for any ball
$\sB_R$.

Next we show that $\sEmin=\Lambda$.
It is obvious that  $\sEmin\le\Lambda$.
Suppose that $\sEmin<\Lambda$.
Let $v\in \bUsm$ (see \cref{T3.1}) and $\sB_r\Supset\cK$.
Then, as shown \cite{ABS19}, the solution
$V$ of \cref{ET3.1A} has the stochastic representation
\begin{equation}\label{ET3.2C}
V(x) \,=\, \Exp_x^v\Bigl[\E^{\int_0^{\uuptau_r}(c_v(X_t)-\sEmin)\, \D{t}}\,
V(X_{\uuptau_r})\Bigr], \quad x\in\sB^c_r\,, \; v\in\bUsm\,.
\end{equation}
Since $\widetilde{V}\in\sorder(\Lyap)$, it follows from \cref{ET3.2B} 
and \cite[Lemma~3.2]{AB18a} that
\begin{equation}\label{ET3.2D}
\widetilde{V}(x)\,\le \, \Exp_x^v\Bigl[\E^{\int_0^{\uuptau_r}(c_v(X_t)-\Lambda)\,\D{t}}\,
\widetilde{V}(X_{\uuptau_r})\Bigr]
\,\le\, \Exp_x^v\Bigl[\E^{\int_0^{\uuptau_r}(c_v(X_t)-\sEmin)\,\D{t}}\,
\widetilde{V}(X_{\uuptau_r})\Bigr], \quad x\in\sB^c_r\,.
\end{equation}
Let $\kappa=\inf_{\sB_r}\frac{V}{\widetilde V}$.
Then $V-\kappa\widetilde V\ge 0$ by \cref{ET3.2C,ET3.2D}, and it vanishes at some
point in $\bar\sB_r$. On the other hand, we have
$$\cA (V-\kappa\widetilde V)\bigl(x, v(x)\bigr) -\sEmin (V-\kappa\widetilde V)
\,\le\, 0\quad \text{in\ } \Rd\,.$$
Hence, by the strong maximum principle we must have $V=\kappa\widetilde{V}$, which
implies that $\sEmin=\Lambda$. Thus we reach a contradiction, and we conclude
that $\sEmin=\Lambda$.
Finally, from \cref{T3.1}\,(c), it follows that $V=\widetilde{V}$.
This also implies that the sequence $\{V_k\}$ converges to $V$.
This completes the proof.
\end{proof}

%%%%%%%%%%%%%%%%%%%%%%%%%%%%%%%%%%%%%%%%%%%%%%%%%%%%%%%%%%%%%%%%%%%%%%%%%%%%%%%%
\section{Policy improvement for the maximization problem}\label{S-maxrisk}

In this section we study the maximization problem, under \cref{A2.2}.
We start the the following theorem.

%%%%%%%%%%%%%%%%%%%%%%%%%%%%%%%%%%%%%%%%%%%%%%%%%%%%%%%%%%%%%%%%%%%%%%%%%%%%%%%%
\begin{theorem}\label{T4.1}
Grant \cref{A2.2}.
There exists a positive solution $\widehat{V}\in\cC^2(\Rd)$ satisfying
\begin{equation}\label{ET4.1A}
\max_{\zeta\in\Act}\, \cA^c \widehat{V}(x,\zeta) \,=\,
\varrho_* \widehat{V}(x)\quad \text{in\ } \Rd\,, \quad \text{and}\quad \widehat{V}(0)=1\,.
\end{equation}
In addition, if $\bUsm\subset\Usm$ denotes the class of Markov controls $v$ which satisfy
$$b\bigl(x,v(x)\bigr)\cdot\grad\widehat{V}(x) + c\bigl(x,v(x)\bigr)\widehat{V}(x)\,=\,
\max_{\zeta\in\Act}\,\bigl[b(x,\zeta)\cdot\grad\widehat{V}(x)
+ c(x,\zeta)\widehat{V}(x)\bigr]
\quad \text{a.e.\ in\ } \Rd\,,$$
then the following hold:
\begin{itemize}
\item[(a)]
$\bUsm= \Usms$, and $\lamstr(c_v)=\sEmax=\varrho_*$ for all $v\in\bUsm$.
Here, $\Usms$ denotes the set of all optimal stationary Markov controls.
\item[(b)]
\Cref{ET4.1A} has a unique solution in $\widehat{V}\in\cC^2(\Rd)$, and
this vanishes at infinity.
\item[(c)] For any $v\in\bUsm$ we have
\begin{equation}\label{ET4.1B}
\widehat{V}(x) \,=\,
\Exp_x^v\left[e^{\int_0^{\uuptau_r} (c_v(X_t)-\varrho_*)\,\D{t}}\,
\widehat{V}(X_{\uuptau_r})\Ind_{\{\uuptau_r<\infty\}}\right] 
\quad \forall\,x\in \bar\sB^c_r\,, \ \forall r>0\,.
\end{equation}
\end{itemize}
\end{theorem}

\begin{proof}
The statement of the theorem is the same as \cite[Theorem~3.1]{ABBK-19}.
However, \cref{A2.2} differs from \cite[Assumption~3.1\,(iii)]{ABBK-19}.
In \cref{EA2.2} if we replace $\sEmax$ with $\varrho_*$, then
the assertions follow from \cite[Theorem~3.1]{ABBK-19}, since any limit of
Neumann eigenvalues cannot be less than the principal eigenvalue.
Suppose then that
\begin{equation}\label{PT4.1A}
\varrho_*\,\le\,\lim_{r\to\infty}\,\max_{(x,\zeta)\in\sB_r^c\times\Act} c(x,\zeta)
\,<\,\sEmax \,<\, \infty\,.
\end{equation}
We use the fact that $\varrho_*$ is a convex, and nondecreasing function of
$c\colon\Rd\times\Act\to\RR$.
Indeed, this follows from the convexity of the maximal operator,
and since the principal eigenvalue $\lamstr(\sL^f)$ is a convex function of $f$.
We add the dependence of $\varrho_*$ on
the coefficient $c$ of the operator explicitly
in the notation by denoting the eigenvalue as $\varrho_*(c)$.
It is also clear from the remark following \cref{A2.2} that
$\varrho_*(c+\delta\Ind_{\sB_1})$ is finite for each $\delta>0$,
and therefore, being convex, it is a continuous function of $\delta\in(0,\infty)$.
Therefore, since $\varrho_*$ is the limit of Dirichlet eigenvalues, which are
strictly increasing as a function of $c$, it follows from
\cref{PT4.1A} that there exists some $\delta>0$ such that 
\begin{equation*}
\lim_{r\to\infty}\,\max_{(x,\zeta)\in\sB_r^c\times\Act} c(x,\zeta)
\,<\, \varrho_*(c+\delta \Ind_{\sB_1}) \,<\, \sEmax\,.
\end{equation*}
But then, as argued earlier, the assertions of the theorem hold,
and these then imply the first equality in 
\begin{equation*}
\sEmax(c+\delta \Ind_{\sB_1})\,=\,\varrho_*(c+\delta \Ind_{\sB_1})\,<\,\sEmax\,.
\end{equation*}
Thus we are led to a contradiction.
This precludes \cref{PT4.1A} as a possibility, and completes the proof.
\end{proof}

Next we  state the PIA.

%%%%%%%%%%%%%%%%%%%%%%%%%%%%%%%%%%%%%%%%%%%%%%%%%%%%%%%%%%%%%%%%%%%%%%%%%%%%%%%%
\begin{algorithm} Policy iteration.
\begin{itemize}
\item[1.]
Initialization. Set $k=0$ and select any $v_0\in\Usm$ which satisfies
$$\lamstr(c_{v_0}) \,>\, 
\lim_{r\to\infty}\,\max_{(x,\zeta)\in\sB_r^c\times\Act} c(x,\zeta)\,.$$
\item[2.]
Value determination. Let $\widehat{V}_k\in \Sobl^{2,p}(\Rd)$, $p>d$,
be the unique principal eigenfunction satisfying
\begin{equation}\label{E4.3}
\trace\bigl(a(x)\grad^2 \widehat{V}_k(x)\bigr) + b(x,v_k)\cdot\grad \widehat{V}_k(x)
+ c(x,v_k) \widehat{V}_k(x)\,=\,
 \lamstr(c_{v_k}) \widehat{V}_k(x) \quad \text{in\ } \Rd\,, \quad \widehat{V}_k(0)=1\,.
\end{equation} 
Existence of unique eigenfunction in \cref{E4.3} follows from \cite{ABBK-19}.
Let 
$\widehat\lambda_k=\lamstr(c_{v_k})$.
\smallskip
\item[3.]
Policy improvement. Choose any $v_{k+1}\in\Usm$ satisfying
$$v_{k+1}(x)\,\in\, \Argmax_{\zeta\in\Act} \left[b(x,\zeta)\cdot\grad \widehat{V}_k(x)
+ c(x,\zeta) \widehat{V}_k(x)\right], \quad x\in\Rd\,.$$
\end{itemize}
\end{algorithm}

The main result of this section is the following.

%%%%%%%%%%%%%%%%%%%%%%%%%%%%%%%%%%%%%%%%%%%%%%%%%%%%%%%%%%%%%%%%%%%%%%%%%%%%%%%%
\begin{theorem}\label{T4.2}
Grant \cref{A2.2}. Then, the following hold:
\begin{itemize}
\item[(i)]
For all $k\in\NN$, we have $\widehat\lambda_k\ge \widehat\lambda_{k-1}$,
and $\lim_{k\to\infty}\widehat\lambda_k=\sEmax$.
\item[(ii)]
The sequence $\{\widehat{V}_k\}$ converges weakly in $\Sobl^{2,p}(\Rd)$,
$p>d$, to the unique solution $\widehat{V}$ of \cref{ET4.1A}.
\end{itemize}
\end{theorem}

We divide the proof in several lemmas.
We adopt the notation in \cref{E-notaA}.

%%%%%%%%%%%%%%%%%%%%%%%%%%%%%%%%%%%%%%%%%%%%%%%%%%%%%%%%%%%%%%%%%%%%%%%%%%%%%%%%
\begin{lemma}\label{L4.1}
We have $\widehat{\lambda}_{k+1}\ge \widehat{\lambda}_{k}$.
\end{lemma}

\begin{proof}
We assume that $\widehat{\lambda}_k> \lim_{r\to\infty}\,\max_{(x,\zeta)\in\sB_r^c\times\Act} c(x,\zeta)$, and establish that
$\widehat{\lambda}_{k+1}\ge \widehat{\lambda}_{k}$.
We employ the method of induction.
This holds for $k=0$.
First, we show that 
\begin{equation}\label{EL4.1A}
\sE_x(c_{k+1}, v_{k+1}) \,\ge\, \widehat{\lambda}_{k} \quad \forall\, x\in\Rd\,.
\end{equation}
Note that
$$\trace\bigl(a(x)\grad^2 \widehat{V}_k(x)\bigr) + b_k(x)\cdot\grad \widehat{V}_k(x)
+ c_k(x)\widehat{V}_k(x)\,=\,\widehat\lambda_k \widehat{V}_k(x)\,,$$
and $\lim_{\abs{x}\to\infty} \widehat{V}_k(x)=0$ by \cref{T4.1}
(see also \cite[Theorem~3.2]{ABBK-19}).
In particular, $\widehat{V}_k\in\cC_b(\Rd)$.
Moreover, 
\begin{equation}\label{EL4.2B}
\begin{aligned}
\trace\bigl(a(x)\grad^2 \widehat{V}_k(x)\bigr)
+ b_{k+1}(x)\cdot&\grad \widehat{V}_k(x)
+ \bigl(c_{k+1}(x)-\widehat{\lambda}_k\bigr) \widehat{V}_k(x)\\
&=\, \max_{\zeta\in\Act}\, \cA^c \widehat{V}_k(x,\zeta)
- \widehat{\lambda}_k \widehat{V}_k(x)\\
&\ge\, \trace\bigl(a(x)\grad^2 \widehat{V}_k(x)\bigr)
+ b_{k}(x)\cdot\grad \widehat{V}_k(x)
+ (c_{k}-\widehat{\lambda}_k) \widehat{V}_k(x)\\
&=0\,.
\end{aligned}
\end{equation}
Recall the definition $\uptau_n= \uptau(\sB_n)$ in \cref{Snot}.
Applying the It\^{o}--Krylov formula to \cref{EL4.2B}, we find that
\begin{equation*}
\Exp_x^{v_{k+1}}
\Bigl[\E^{\int_0^{T\wedge\uptau_n}(c_{k+1}(X_t)-\widehat{\lambda}_k)\,\D{t}}\,
\widehat{V}_{k}(X_{T\wedge\uptau_n})\Bigr]\,\ge\, \widehat{V}_k(x)
\qquad\forall\, x\in\sB_n\,.
\end{equation*}
Since $\norm{c^+}_\infty<\infty$, letting $n\to\infty$ above,
and applying the dominated convergence theorem, we obtain
$$\Exp_x^{v_{k+1}}\Bigl[\E^{\int_0^{T}(c_{k+1}(X_t)-\widehat{\lambda}_k)\,\D{t}}\,
 \widehat{V}_{k}(X_{T})\Bigr]\,\ge\, \widehat{V}_k(x)\,.$$
Thus\,
\begin{align*}
\log \widehat{V}_k(x)\,\le\, -\widehat{\lambda}_k T + \log\Exp_x^{v_{k+1}}
\Bigl[\E^{\int_0^{T}c_{k+1}(X_t)\D{t}}\Bigr]
+ \log \norm{\widehat{V}_k}_\infty\,.
\end{align*}
Now dividing by $T$ on both sides and letting $T\to\infty$ we have \cref{EL4.1A}.

To complete the proof, we show that
$\widehat{\lambda}_{k+1}= \sup_{x\in\Rd}\sE_x(c_{k+1}, v_{k+1})$.
In view of \cref{EL4.2B}, the calculations above,
and \cite[Theorem~1.7]{Berestycki-15}, we note that
$\widehat{\lambda}_{k+1}\le \sup_{x\in\Rd}\sE(c_{k+1}, v_{k+1})$.
Arguing as in the proof of \cref{T4.1}, if
\begin{equation}\label{EL4.2C}
\widehat{\lambda}_{k+1}\,>\,
\lim_{r\to\infty}\,\max_{(x,\zeta)\in\sB_r^c\times\Act} c(x,\zeta)\,,
\end{equation}
then we have 
$\widehat{\lambda}_{k+1} = \sE_x(c_{k+1}, v_{k+1})$ for all $x$.
Now suppose that 
$$\widehat{\lambda}_{k+1} \,\le\, 
\lim_{r\to\infty}\,\max_{(x,\zeta)\in\sB_r^c\times\Act} c(x,\zeta)
 \,<\, \sup_{x\in\Rd}\sE(c_{k+1}, v_{k+1})\,,$$
where the last inequality follows from
\cref{EL4.1A}. We know from \cite{Berestycki-15} that
$\delta\mapsto \lamstr(c_{k+1}+\delta\Ind_{\sB_1})$ is a convex function.
Again, since $\lamstr(c_{k+1}+\delta\Ind_{\sB_1})$ are obtained as a increasing
limit of Dirichlet principal eigenvalues, it follows that 
$$\lim_{\delta\to\infty} \lamstr(c_{k+1}+\delta\Ind_{\sB_1})
\,\ge\, \lim_{\delta\to\infty} \lambda_1(c_{k+1}+\delta) \,=\, \infty\,,$$
where $\lambda_1(c_{k+1}+\delta)$ denotes the Dirichlet principal eigenvalue
in the unit ball.
Thus, we can find a $\delta_\circ$ satisfying
$$\lim_{r\to\infty}\,\max_{(x,\zeta)\in\sB_r^c\times\Act} c(x,\zeta)
\,<\,\lamstr(c_{k+1}
+\delta_\circ\Ind_{\sB_1})\,<\,\sup_{x\in\Rd}\sE(c_{k+1}, v_{k+1})\,.$$
Therefore, as argued earlier, it follows that 
$$\sup_{x\in\Rd}\sE(c_{k+1}+\delta_\circ\Ind_{\sB_1}, v_{k+1}) \,=\, \lamstr(c_{k+1}
+\delta_\circ\Ind_{\sB_1}) \,<\, \sup_{x\in\Rd}\sE(c_{k+1}, v_{k+1})\,,$$
which is a contradiction.
Thus, \cref{EL4.2C} must hold,
and this completes the proof.
\end{proof}

Next we establish the strict monotonicity of the
eigenvalue  at every $c_{k}$ (see \cref{D2.3}).
Recall the definition of $\sL_k$ in \cref{E-notaA}.

%%%%%%%%%%%%%%%%%%%%%%%%%%%%%%%%%%%%%%%%%%%%%%%%%%%%%%%%%%%%%%%%%%%%%%%%%%%%%%%%
\begin{lemma}\label{L4.2}
For every $k\ge 0$, the principal eigenvalue of $\sL_k$ is strictly
monotone with respect to the potential $c_k$.
\end{lemma}

\begin{proof}
Let $h\neq 0$ be a nonnegative, compactly supported continuous function.
From the definition
of the principal eigenvalue it is clear that
$\lambda_1\df\lamstr(\sL_k- h)\le \lamstr(\sL_k)\df\lambda_2$.
Suppose that $\lambda_1=\lambda_2$.
Let $\Psi_i\in \Sobl^{2,p}(\Rd)$ be the eigenfunction corresponding to the eigenvalue
$\lambda_i$ for $i=1,2$.
Now, using \cref{L4.1} it follows that
$$\lambda_1 \,=\, \lambda_2 \,>\,
\lim_{r\to\infty}\,\max_{(x,\zeta)\in\sB_r^c\times\Act}\,c(x,\zeta)\,,$$
and therefore, $\lim_{\abs{x}\to\infty}\Psi_i(x)=0$ and both the eigenfunctions have 
stochastic representation \cref{ET4.1B}.
Since $h$ is compactly supported, by choosing $r$ large enough we obtain
\begin{align*}
\Psi_i(x) \,=\, \Exp_x^{v_k}\Bigl[e^{\int_0^{\uuptau_r} (c_k(X_t)-\lambda_1)\D{t}}\,
\Psi_i(X_{\uuptau_r})\Ind_{\{\uuptau_r<\infty\}}\Bigr], 
\quad x\in \bar\sB^c_r\,,\ i=1,2\,.
\end{align*}
Now choose a constant $\kappa>0$ so that $\kappa\Psi_2\ge \Psi_1$,
and equality holds at some
point in $\bar\sB_r$. On the other hand, we also have
$$\sL_k(\kappa\Psi_2-\Psi_1) \,=\, -h\Psi_1\le 0\,.$$
Thus, by the strong maximum principle we must have $\kappa \Psi_2=\Psi_1$ in $\Rd$,
which contradicts the fact that $h\ne 0$.
Therefore, $\lambda_1<\lambda_2$, thus completing the proof.
\end{proof}

Recall that the twisted process is given by
\begin{equation}\label{E4.7}
\D{Y^k_t} \,=\, b_k(Y^k_t)\,\D{t} + 2a(Y^k_t) \grad \log\widehat{V}_k(Y^k_t)\,\D{t}
+ \upsigma(Y^k_t)\,\D{W_t}\,.
\end{equation}
From \cref{T2.1}, \cref{R2.2,L4.2} we note that \cref{E4.7} is exponentially ergodic.
Let $\Tilde\upmu_k$ be the corresponding unique invariant measure.
In the lemma which follows, we show that $\{\Tilde\upmu_k\}_{k\ge 0}$ is tight.

%%%%%%%%%%%%%%%%%%%%%%%%%%%%%%%%%%%%%%%%%%%%%%%%%%%%%%%%%%%%%%%%%%%%%%%%%%%%%%%%
\begin{lemma}\label{L4.3}
The family of invariant measures $\{\Tilde\upmu_k\}_{k\ge 0}$ is tight.
\end{lemma}

\begin{proof}
Let $\epsilon, r>0$ be such that
$\max_{\zeta\in\Act}\,c(x,\zeta)-\widehat\lambda_0<-\epsilon$ for all $x\in \sB^c_r$.
It follows by \cref{L4.1}, that
$\max_{\zeta\in\Act}\,c(x,\zeta)-\widehat\lambda_k<-\epsilon$
for all $x\in \sB^c_r$ and $k\ge 0$.
Let $\widetilde\sL_k$ denote the extended generator of \cref{E4.7}.
Then an easy calculation reveals that, with $\Check{V}_k\df(\widehat{V}_k)^{-1}$,
we have
\begin{equation}\label{EL4.3A}
\widetilde\sL_k \check{V}_k + (\widehat\lambda_k-c_k)\check{V}_k \,=\, 0
\quad \text{in\ } \Rd\,,
\end{equation}
and moreover, $\check{V}_k$ is inf-compact.
As done earlier,
denote by $\widetilde\Exp^k$ the expectation operator
on the canonical space of the process \cref{E4.7}.
Then applying the It\^{o}--Krylov formula and Fatou's lemma, we obtain
\begin{align*}
\check{V}_k(x) &\,\ge\, \widetilde\Exp_x^k\Bigl[\E^{\int_0^{\uuptau_r}
(\widehat\lambda_k-c_{k}(Y^k_t))\,\D{t}}\, \check{V}_k(Y^k_{\uuptau_r})\Bigr]\\
&\,\ge\,\Bigl(\max_{\bar\sB_r} \widehat{V}_k\Bigr)^{-1}\,
\widetilde\Exp_x^k\bigl[\E^{\epsilon \uuptau_r}  \bigr]\\
&\,\ge\,\Bigl(\max_{\bar\sB_r} \widehat{V}_k\Bigr)^{-1}\, \epsilon\,
\exp \bigl(\widetilde\Exp_x^k[ \uuptau_r]\bigr)\qquad\forall\,x\in\bar\sB^c_r\,.
\end{align*}
Since, by \cref{L3.2}, $\abs{\grad \log\widehat{V}_k}\le \kappa(1+ \abs{x})$
for all $x\in\Rd$ and $k\in\NN$, it follows from 
\hyperlink{A2}{(A2)} that $\check{V}_k$ is inf-compact,
uniformly in $k$ (see for instance, \cite[Lemma~5.1]{AB-19}),
that is, $\inf_k \check{V}_k$ is inf-compact.
Then, the result follows from \cref{EL4.3A}.
\end{proof}

Next, we present the proof of \cref{T4.2}.

%%%%%%%%%%%%%%%%%%%%%%%%%%%%%%%%%%%%%%%%%%%%%%%%%%%%%%%%%%%%%%%%%%%%%%%%%%%%%%%%
\begin{proof}[Proof of \cref{T4.2}]
Let $\widehat\Lambda=\lim_{k\to\infty}\widehat{\lambda}_k$.
Existence follows from \cref{L4.1}.
In view of the proof of \cref{T3.2}, it is enough to 
show that $\widehat\Lambda=\sEmax$.
As earlier, we define
$$\max_{\zeta\in\Act}\,\cA^c \widehat{V}_k(x,\zeta)
-\widehat\lambda_k \widehat{V}_k(x)\,=\,\widehat\psi_{k+1}(x) \widehat{V}_k(x)\,.$$
Then $\widehat{\psi}_k$ satisfies \cref{E3.8}, and using \cref{L4.3}, we obtain
$$\lim_{k\to\infty}\int_{\sB_R}\widehat{\psi}_k(x)\, \D{x}\,=\,0
\quad \forall\, R>0\,.$$
Also, the uniform estimate on $\check{V}_k$ (see the proof of \cref{L4.3})
shows that any limit of $\{\widehat{V}_k\}$ must vanish at infinity.
Therefore, we can follow the arguments in \cref{T3.2} together with \cref{T4.1}
to complete the proof.
\end{proof}

%%%%%%%%%%%%%%%%%%%%%%%%%%%%%%%%%%%%%%%%%%%%%%%%%%%%%%%%%%%%%%%%%%%%%%%%%%%%%%%%
\section{A general  result on convergence}\label{S5}

In this section we provide sufficient conditions for the PIA to converge,
without assuming blanket stability hypotheses or near-monotonicity
of the running cost.
We apply these to the minimization problem under a near-monotone cost
hypothesis in \cref{S5.1}.
We address the minimization problem.
Let
\begin{equation*}
\cG f(x) \,\df\, \trace\bigl(a(x)\grad^2 f(x)\bigr)
+ \min_{\zeta\in\Act}\,\bigl[ b(x,\zeta)\cdot\grad{f}(x) + c(x,\zeta){f}(x)\bigr]\,,
\end{equation*}
and denote by $\lamstr(\cG)$ the generalized principal eigenvalue of the operator
$\cG$ on $\Rd$, which is defined by \cref{E-lamstr}
with $\sL^f$ replaced by $\cG$, and which is assumed to be finite.
In this section, the coefficients $a$, $b$, and $c$ are not restricted to satisfy
\hyperlink{A1}{(A1)}--\hyperlink{A4}{(A4)}.
Rather, we assume that they satisfy the more general
hypotheses in \cite[Section~2.1]{AB18a}.
That is, we replace \hyperlink{A1}{(A1)}--\hyperlink{A4}{(A4)}
with the following:

\begin{itemize}
\item[\hypertarget{B1}{(B1)}]
The functions
$b\colon\RR^{d}\times\Act\to\RR^{d}$ and 
$\upsigma\colon\RR^{d}\to\RR^{d\times d}$
are continuous and satisfy
\begin{equation*}
\abs{b(x,u)-b(y, u)} + \norm{\upsigma(x) - \upsigma(y)}
\,\le\,C_{R}\,\abs{x-y}\qquad\forall\,x,y\in B_R\,,\ \forall\, u\in\Act\,,
\end{equation*}
for some constant $C_{R}>0$ depending on $R>0$.

\item[\hypertarget{B2}{(B2)}]
For some $C_0>0$, we have
\begin{equation*}
\sup_{\zeta\in\Act}\, \langle b(x,\zeta),x\rangle^{+} + \norm{\upsigma(x)}^{2}\,\le\,C_0
\bigl(1 + \abs{x}^{2}\bigr) \qquad \forall\, x\in\RR^{d}\,.
\end{equation*}

\item[\hypertarget{B3}{(B3)}]
For each $R>0$, it holds that
\begin{equation*}
\sum_{i, j=1}^d a^{ij}(x) \xi_i \xi_j \ge C^{-1}_R \abs{\xi}^2 \qquad \forall x\in B_R\,,
\end{equation*}
and for all $\xi=(\xi_1, \ldots, \xi_d)\transp\in\Rd$,
with $C_R$ the constant in \hyperlink{B1}{(B1)}.

\item[\hypertarget{B4}{(B4)}]
The running cost $c(x,\zeta)$ is in $\cC(\Rd\times\Act, \RR)$, is bounded below in $\Rd$,
and is locally
Lipschitz in $x$ uniformly with respect to $\zeta\in\Act$.
\end{itemize}

Hypotheses \hyperlink{B1}{(B1)}--\hyperlink{B4}{(B4)} are assumed throughout
this section.

\begin{definition}
We adopt the notation
\begin{equation*}
\cA^c_v f(x)\,\df\, \trace\bigl(a(x)\grad^2 f(x)\bigr)
+ b_v(x)\cdot\grad f(x) + c_v(x) f(x)\,,\quad
v\in\Usm\,,
\end{equation*}
where $b_v$ and $c_v$ are as in \eqref{Esimpl},
and let $\lamstr_v$ denote the principal eigenvalue $\lamstr(\cA^c_v)$ of
$\cA^c_v$ on $\Rd$ defined in \eqref{E-lamstr}.
We let
\begin{equation*}
\Lambda^* \,\df\, \inf_{v\in\Usm}\, \lamstr_v\,,
\end{equation*}
and
\begin{equation*}
\mathcal{U}^* \,\df\, \bigl\{v\in\Usm\colon \lamstr_v=\Lambda^*\bigr\}\,.
\end{equation*}
\end{definition}

We denote $\lamstr(\cG)$ also as $\lamstr(\cG,c)$
when we need its dependence on the coefficient
$c$ of the operator to be explicitly captured in the notation.
The following theorem extends the results in \cite{AB18}.

%%%%%%%%%%%%%%%%%%%%%%%%%%%%%%%%%%%%%%%%%%%%%%%%%%%%%%%%%%%%%%%%%%%%%%%%%%%%%%%%
\begin{theorem}\label{T5.1}
Suppose that the eigenvalue $\lamstr(\cG,c)$ of the semilinear
operator $\cG$ is strictly monotone at $c$ on the right, that is,
for any non-trivial
nonnegative function $h\colon\Rd\to\RR$, we have
$\lamstr(\cG,c)<\lamstr(\cG,c+h)$ (compare with \cref{D2.3}).
Then there exists a unique positive $\Psi_*\in\cC^2(\Rd)$, which satisfies
$\Psi_*(0)=1$, and
\begin{equation}\label{ET5.1A}
\trace\bigl(a(x)\grad^2 \Psi_*(x)\bigr)
+ \min_{\zeta\in\Act}\,\bigl[b(x,\zeta)\cdot\grad{\Psi_*}(x) + c(x,\zeta){\Psi_*}(x)\bigr]
\,=\,\lamstr(\cG) \Psi_*(x)\,.
\end{equation}
In addition, if $\bUsm\subset\Usm$ denotes
the class of Markov controls $v$ which satisfy
\begin{equation*}
b_v(x)\cdot\grad\Psi_*(x) + c_v(x)\Psi_*(x)\,=\,
\min_{\zeta\in\Act}\,\bigl[ b(x,\zeta)\cdot\grad\Psi_*(x) + c(x,\zeta)\Psi_*(x)\bigr]
\quad \text{a.e.\ in}\ \Rd\,,
\end{equation*}
then $\bUsm=\mathcal{U}^*$.
\end{theorem}

\begin{proof}
First note that, under the hypothesis
that $\lamstr(\cG)$ is finite, existence of a positive solution
$\Psi_*$ is standard,
and can be constructed as a limit of Dirichlet eigenvalue problems
on an increasing sequence of balls in $\Rd$ as in \cite[Lemma~3.1]{AB18}.
Indeed, by \cite[Lemma~3.1]{AB18},
there exists a unique pair
$(\widehat{\Psi}_{n},\Hat\lambda_n)
\in\bigl(C^2(\sB_{n})\cap\cC(\Bar{\sB}_{n})\bigr)\times\RR$, $n\in\NN$, satisfying
$\widehat{\Psi}_{n}>0$ on $\sB_{n}$,
$\widehat{\Psi}_{n}=0$ on $\partial \sB_{n}$, and $\widehat{\Psi}_{n}(0)=1$,
which solves
\begin{equation*}
\min_{\zeta\in\Act}\,
\bigl[\cA \widehat{\Psi}_{n}(x,\zeta) + c(x,\zeta)\,\widehat{\Psi}_{n}(x)\bigr]
\,=\, \Hat\lambda_n\,\widehat{\Psi}_{n}(x)\,,\qquad x\in \sB_{n}\,,
\end{equation*}
and $\Hat\lambda_n$ is an increasing sequence.
Any limit point $\widehat{\Psi}_{n}$ as $n\to\infty$
satisfies \cref{ET5.1A}.
However, here we use a different construction for $\Psi_*$.
Since
$\Hat\lambda_n<\lamstr(\cG)$ for all $n\in\NN$,
for any $\alpha_n>0$, the Dirichlet problem
\begin{equation}\label{PT5.1B}
\min_{\zeta\in\Act}\,
\bigl[\cA\varphi_{n}(x,\zeta)+ \bigl(c(x,\zeta)-\lamstr(\cG)\bigr)\,\varphi_{n}(x)\bigr]
\,=\, -\alpha_n\,\Ind_{\sB}(x)\quad\text{a.e.\ }x\in B_n\,,\quad
\varphi_{n}=0\text{\ \ on\ \ }\partial B_n\,,
\end{equation}
has a unique positive solution $\varphi_{n}\in\Sobl^{2,p}(B_n)\cap\cC(\Bar{B}_n)$,
for any $p\ge1$, by \cite[Theorem~2.3\,(ii)]{Armstrong-09}
(see also \cite[Theorem~1.1\,(ii)]{Yoshimura-06}). 
It is also easy to see 
that $\varphi_n(0)$ is a continuous increasing function of the constant
$\alpha_n$.
We choose the constant $\alpha_n$ as follows:
if $\Tilde\alpha_n>0$ is such that the solution
$\varphi_n$ of \eqref{PT5.1B} with $\alpha_n=\Tilde\alpha_n$
satisfies $\varphi_n(0)=1$, we set $\alpha_n=\min(1,\Tilde\alpha_n)$.
Passing to the limit as $n\to\infty$ along some subsequence in \cref{PT5.1B},
we obtain a positive solution $\Psi_*\in\Sobl^{2,p}(\Rd)$, for $p>d$, of
\begin{equation}\label{PT5.1C}
\min_{\zeta\in\Act}\,
\bigl[\cA\Psi_*(x,\zeta)+ \bigl(c(x,\zeta)-\lamstr(\cG)\bigr)\,\Psi_*(x)\bigr]
\,=\, -\alpha\,\Ind_{\sB}(x)\,,\qquad x\in\Rd\,,
\end{equation}
for some $\alpha\ge0$.
Writing \cref{PT5.1C} as 
\begin{equation*}
\min_{\zeta\in\Act}\,
\Bigl[\cA\Psi_*(x,\zeta)+ \Bigl(c(x,\zeta)
+\alpha\tfrac{\Ind_{\sB}(x)}{\Psi_*(x)}-\lamstr(\cG)\Bigr)\,\Psi_*(x)\Bigr]
\,=\, 0\,,\qquad x\in\Rd\,,
\end{equation*}
it is clear that the strict right-monotonicity of $\lamstr(\cG,c)$
at $c$ implies that $\alpha=0$.
Thus $\Psi_*$ solves \cref{ET5.1A}, and is in $\cC^2(\Rd)$
by elliptic regularity \cite[Theorem~9.19]{GilTru}.

If $\Bar{v}\in\bUsm$, then we have
\begin{equation}\label{PT5.1D}
\cA^c_{\Bar{v}} \Psi_* \,=\, \lamstr(\cG) \Psi_*\quad\text{a.e.\ on\ }\Rd\,.
\end{equation}
It follows from the definition of the principal eigenvalue $\lamstr(\cG)$
that $\lamstr_v=\lamstr(\cA^c_{v})\ge\lamstr(\cG)$ for any $v\in\Usm$.
Therefore, $\lamstr_{\Bar{v}} = \lamstr(\cG)$ for any $\Bar{v}\in\bUsm$
by \cref{PT5.1D}.
This shows that $\Lambda^*=\lamstr(\cG)$ and that $\bUsm\subset\mathcal{U}^*$.

Let $v\in\mathcal{U}^*$.
By It\^o's formula applied to \eqref{PT5.1C}, with $\uuptau_r=\uptau(\sB^c_r)$
as defined in \cref{Snot} for $r<n$,
we obtain
\begin{multline*}
\varphi_n(x) \,\le\, \Exp_x^v\Bigl[\E^{\int_{0}^{\uuptau_r}
[c_v(X_s)-\Lambda^*]\,\D{s}}\,
\varphi_n(X_{\uuptau_r})\,\Ind_{\{\uuptau_r<T\wedge\uptau_n\}}\Bigr]\\[5pt]
+\Exp_x^v\Bigl[\E^{\int_{0}^{T}[c_v(X_s)-\Lambda^*]\,\D{s}}\,
\varphi_n(X_{T})\,
\Ind_{\{T<\uuptau_r\wedge\uptau_n\}}\Bigr]
\qquad\forall\,x\in \sB_n\setminus\sB_r\,,\ \forall\,T>0\,.
\end{multline*}
Then we use the argument in the proof of \cite[Lemma~2.11]{AB18},
by replacing $f$ with $c_v$, and $\Lambda(f)$ with
$\Lambda^*$ in \cite[(2.34)--(2.36)]{AB18}
to obtain
\begin{equation}\label{PT5.1E}
\Psi_*(x) \,\le\, \Exp_x^v
\Bigl[\E^{\int_{0}^{\uuptau_r}[c_v(X_s)-\Lambda^*]\,\D{s}}\,
\Psi_*(X_{\uuptau_r})\,\,\Ind_{\{\uuptau_r<\infty\}}\Bigr]
\qquad\forall\,x\in \sB^c_r\,,
\end{equation}
and for all $v\in\mathcal{U}^*$.
We alert the reader to the fact that this part of the argument
in the proof of \cite[Lemma~2.11]{AB18}
does not rely on the near-monotonicity of $c$.

Next, suppose $v^*\in\mathcal{U}^*$, and
let $\Phi_*\in\Sobl^{2,p}(\Rd)$, $p>d$, be an eigenfunction
for the operator $\cA^c_{v^*}$ corresponding to the
eigenvalue $\lamstr_{v^*}=\Lambda^*$, and such that $\Phi_*(0)=1$.
It also follows from \cref{D2.3} that the strict monotonicity
of $\lamstr(\cG,c)$ at $c$ on the right implies the same property
for the eigenvalue $\lamstr_{v^*}$ of the linear operator $\cA^c_{v^*}$
for any $v^*\in\mathcal{U}^*$.
This implies that $\lamstr_{v^*}$ is a simple eigenvalue
of the linear operator $\cA^c_{v^*}$ for
$v^*\in\mathcal{U}^*$, and 
that $\Phi_*$ has the stochastic representation
\begin{equation}\label{PT5.1F}
\Phi_*(x)\;=\;\Exp_x^{v^*}\Bigl[\E^{\int_0^{\uuptau_r}[c_{v^*}(X_s)-\Lambda^*]\, \D{s}}\,
\Phi_*(X_{\uuptau})\,\Ind_{\{\uuptau_r<\infty\}}\Bigr]\qquad
\forall\, x\in \Bar{\sB}_r^c\,,\quad\forall\,r>0
\end{equation}
by \cite[Theorem~2.1]{ABS19}.
A standard application of the strong maximum principle as in the
proof of \cite[Lemma~3.5]{AB18} using \cref{PT5.1E,PT5.1F} then shows
that $\Phi_*=\Psi_*$.
This shows that $\mathcal{U}^*\subset\bUsm$, and thus we must have equality.
Uniqueness of the solution $\Psi_*$ of \cref{ET5.1A} clearly then follows
from the simplicity of $\lamstr_{v^*}$ for
any $v^*\in\mathcal{U}^*$.
This completes the proof.
\end{proof}

Let $\cM_*$ denote the class of invariant probability measures, corresponding
to the ground state processes with $v\in\bUsm$.
More precisely, $\cM_*$ is the collection of the invariant probability
measures corresponding to the extended generators
\begin{equation*}
\widetilde\sL_v f(x) \,\df\, \trace\bigl(a(x)\grad^2 f(x)\bigr)
+ \bigl(b\bigl(x,v(x)\bigr) + 2a(x)\grad\log\Psi_*(x)\bigr) \cdot\grad f(x)\,,
\quad v\in\bUsm\,.
\end{equation*}

Note that for $\cM_*\ne\varnothing$, it is sufficient that
the eigenvalue $\lamstr_v$ of the operator $\cA^c_v$
be strictly monotone at $c_v$ 
for some $v\in\bUsm$ \cite[Theorem~2.2]{ABS19}.

The following assumption is enforced throughout this section, without
further mention.

%%%%%%%%%%%%%%%%%%%%%%%%%%%%%%%%%%%%%%%%%%%%%%%%%%%%%%%%%%%%%%%%%%%%%%%%%%%%%%%%
\begin{assumption}\label{A5.1}
The eigenvalue $\lamstr(\cG)$ of the semilinear
operator $\cG$ is strictly monotone at $c$ on the right,
and $\cM_*\ne\varnothing$.
\end{assumption}

We let $\bm\Psi_v$ denote the set of eigenfunctions $\Psi_v$ obtained
as limits of Dirichlet eigenvalue problems
for the operator $\cA^c_v$
, normalized so that $\Psi_v(0)=1$,
and for $\rho\in\RR$,
we define the space of functions
\begin{equation*}
\bm\Psi(\rho) \,\df\,
\bigl\{V\in\bm\Psi_v\,\colon \lamstr_v\le\rho\,,\;v\in\Usm\bigr\}\,,
\end{equation*}
that is, the set of eigenfunctions corresponding to eigenvalues not exceeding $\rho$.
It is clear that $\bm\Psi(\rho)=\varnothing$ if $\rho<\lamstr(\cG)$,
and,  under \cref{A5.1}, $\bm\Psi\bigl(\lamstr(\cG)\bigr)=\{\Psi_*\}$
by \cref{T5.1}.

We let
\begin{equation*}
\mathcal{U}(\rho) \,\df\, \{v\in\Usm\,\colon \lamstr_v\le\rho\}\,.
\end{equation*}

We need the following definition.
%%%%%%%%%%%%%%%%%%%%%%%%%%%%%%%%%%%%%%%%%%%%%%%%%%%%%%%%%%%%%%%%%%%%%%%%%%%%%%%%
\begin{definition}\label{D5.2}
Let $\rho\in(\lamstr(\cG),\infty)$.
\begin{enumerate}
\item[\ttup{i}]
We say that $\rho$ has \hypertarget{PA}{\emph{Property~A}} if
\begin{equation}\label{ED5.2A}
\int_\Rd\biggl(\sup_{V\in\bm\Psi(\rho)}\,\frac{V(x)}{\Psi_*(x)}\biggr)\,
\mu(\D{x})\,<\,\infty
\end{equation}
for some $\mu\in\cM_*$.
\item[\ttup{ii}]
We say that $\rho$ has \hypertarget{PB}{\emph{Property~B}} if
the ground state diffusions with generators $\widetilde\sL_v$,
$v\in\mathcal{U}(\rho)$,
are positive recurrent and the corresponding invariant probability measures
$\bigl\{\Tilde\upmu_v\,\colon v\in\mathcal{U}(\rho)\bigr\}$
are tight.
\end{enumerate}
\end{definition}

\hyperlink{PB}{Property~B} implies that $\bm\Psi_v$ is a singleton
for all $v\in\mathcal{U}(\rho)$
\cite[Lemma~2.7 and Theorem~2.3]{ABS19}.
Recall \cref{Alg3.1}.
We have the following convergence result.

%%%%%%%%%%%%%%%%%%%%%%%%%%%%%%%%%%%%%%%%%%%%%%%%%%%%%%%%%%%%%%%%%%%%%%%%%%%%%%%%
\begin{lemma}\label{L5.1}
Suppose that $\rho\in(\lamstr(\cG),\infty)$ has
\hyperlink{PB}{Property~B}.
Then, if the control $v_0\in\Usm$ in the initialization of \cref{Alg3.1} is
such that $\lamstr_{v_0}\le\rho$, then
 $\lambda_k\to\Hat\lambda\in[\lamstr(\cG),\rho]$, and
the sequence $\{V_k\}$ converges uniformly on compact sets to some
$\Phi\in\cC^2(\Rd)$ satisfying
\begin{equation}\label{EL5.1A}
\min_{\zeta\in\Act}\, \cA^c \Phi(x,\zeta) \,=\, \Hat\lambda \Phi(x)\,.
\end{equation}
\end{lemma}

\begin{proof}
We follow the proof of \cref{L3.4}.
It is clear that $\lambda_k\le\lambda_{k-1}$ for all $k\in\NN$,
so that $\lambda_k$ converges to some $\Hat\lambda\in[\lamstr(\cG),\rho]$.
Using \cref{PL3.4C} in \cref{PL3.4A}, and evaluating
at $x=0$, and $T=\infty$, we have
\begin{equation}\label{PL5.1A}
\int_0^\infty \E^{-\Tilde\lambda_k t}\,\widetilde\Exp_0^{k}\bigl[h_k(Y_t^k)\bigr] \D{t}
\,\le\,1\,,
\end{equation}
and we know that $\Tilde\lambda_k\searrow0$.
Let $R>0$, which is used in the proof as a parameter.
Define
\begin{equation*}
J_k(x) \,\df\,
\int_0^\infty \E^{-\Tilde\lambda_k t}\,\widetilde\Exp_0^{v_k}\bigl[h_k(Y_t)\,
\Ind_{\sB_R}(Y_t)\bigr] \D{t}\,.
\end{equation*}
Let $\Tilde\upmu_k$ denote the invariant measure of the ground state process
$\process{Y^k}$.
These are tight by the hypothesis in \hyperlink{PB}{Property~B}.
Thus, by invariance, the Harnack property of the densities of the invariant
measures together with tightness, and the fact that $h_k$ is bounded on each ball $\sB_R$
uniformly in $k$, there
exists a positive constant $\epsilon_1$ depending only on $R>0$, such that
\begin{equation}\label{PL5.1B}
\int_\Rd \Tilde\lambda_k J_k(x)\,\Tilde\upmu_k(\D{x}) \,=\,
\int_{\sB_R}h_k(x)\,\Tilde\upmu_k(\D{x}) \,\ge\, \epsilon_1 \norm{h_k}_{\Lp^1(\sB_R)}\,.
\end{equation}
It is well known that the supremum of $J_k(x)$ on $\Rd$ is realized
at $\Bar\sB_R$ \cite[Lemma~3.6.1]{book}.
Therefore
\begin{equation}\label{PL5.1C}
\int_\Rd  J_k(x)\,\Tilde\upmu_k(\D{x}) \,\le\,
\sup_{\sB_R}\,J_k\,.
\end{equation}
As shown in \cite[(3.7.24)]{book}, by 
employing \cite{AGM99}*{Corollary~2.2}, the resolvent $J_k$ has the
Harnack property.
Thus, there exists a constant $C_{\mathsf{H}}$, independent of $k$, such that
\begin{equation}\label{PL5.1D}
C_{\mathsf{H}}\, \inf_{\sB_R}\,J_k\,\ge\,
\sup_{\sB_R}\,J_k\,.
\end{equation}
Combining \cref{PL5.1B,PL5.1C,PL5.1D}, we have
\begin{equation}\label{PL5.1E}
\begin{aligned}
\norm{h_k}_{\Lp^1(\sB_R)} &\,\le\, \frac{\Tilde\lambda_k}{\epsilon_1}\,\sup_{\sB_R}\,J_k\\
&\,\le\, \frac{\Tilde\lambda_k}{\epsilon_1}\,C_{\mathsf{H}}\, \inf_{\sB_R}\,J_k\\
&\,\le\, \frac{\Tilde\lambda_k}{\epsilon_1}\,C_{\mathsf{H}}\, J_k(0)\\
&\,\le\, \frac{\Tilde\lambda_k}{\epsilon_1}\,C_{\mathsf{H}}\,,
\end{aligned}
\end{equation}
where in the last inequality we use \cref{PL5.1A}.
Since $\Tilde\lambda_k\to0$ as $k\to\infty$,
it follows from \cref{PL5.1E} that $h_k$ converges to $0$ in $\Lp^1$
on every ball $\sB_R$,
and thus also converges in $\Lp^p$ for any $p\ge1$, since it is uniformly bounded
on each ball.
The rest follows exactly as in \cref{L3.4}.
\end{proof}

%%%%%%%%%%%%%%%%%%%%%%%%%%%%%%%%%%%%%%%%%%%%%%%%%%%%%%%%%%%%%%%%%%%%%%%%%%%%%%%%
\begin{lemma}\label{L5.2}
In addition to the hypotheses of \cref{L5.1}, suppose that
$\rho\in(\lamstr(\cG),\infty)$ has \hyperlink{PA}{Property~A}.
Then, $(\Hat\lambda,\Phi)=(\lamstr(\cG),\Psi_*)$ is the unique solution of
\cref{EL5.1A}.
\end{lemma}

\begin{proof}
Suppose that $\Phi\in\Sobl^{2,p}(\Rd)$, $p>d$, solves \cref{EL5.1A} for
some $\Hat\lambda\in[\lamstr(\cG),\rho]$.
Since $\Phi$ is a limit of a sequence $V_k$ of \cref{Alg3.1}, it is
clear that $\Phi\le \sup_{V\in\bm\Psi(\rho)}\,V$.
Let $v\in\bUsm$ be such that the corresponding ground process $\process{Y}$
has an invariant probability measure $\Tilde\upmu_{v}\in\cM_*$, satisfying \cref{ED5.2A}.
We have
\begin{equation*}
\widetilde\sL_v \biggl(\frac{\Phi}{\Psi_*}\biggr)
\,\ge\, \bigl(\Hat\lambda-\lamstr(\cG)\bigr) \biggl(\frac{\Phi}{\Psi_*}\biggr)\,.
\end{equation*}
Therefore, $\bigl\{\frac{\Phi}{\Psi_*}(Y_t)\bigr\}_{t\ge0}$
is a nonnegative submartingale,
and is integrable under the invariant probability measure
$\Tilde\upmu_{v}$.
Then since $\widetilde\Exp_x^{v}\bigl[\tfrac{\Phi}{\Psi_*}(Y_t)\bigr]
\ge \tfrac{\Phi}{\Psi_*}(x)$, we obtain
$\tfrac{\Phi}{\Psi_*}(x)\le \Tilde\upmu_{v}\bigl(\tfrac{\Phi}{\Psi_*}\bigr)$,
so that $\tfrac{\Phi}{\Psi_*}$ is bounded.
Thus $\bigl\{\frac{\Phi}{\Psi_*}(Y_t)\bigr\}_{t\ge0}$ converges a.s.,
and since the process is recurrent it must converge to a constant.
This shows that $\Phi=\kappa\Psi_*$ for some $\kappa>0$.
In turn, this implies that $\Hat\lambda=\lamstr(\cG)$.
\end{proof}

We are ready to state the main convergence result.

%%%%%%%%%%%%%%%%%%%%%%%%%%%%%%%%%%%%%%%%%%%%%%%%%%%%%%%%%%%%%%%%%%%%%%%%%%%%%%%%
\begin{theorem}\label{T5.2}
Grant \cref{A5.1}, and suppose that
$\rho\in(\lamstr(\cG),\infty)$ has Properties~A and B.
Then, provided that the control $v_0\in\Usm$ in the initialization of \cref{Alg3.1} is
such that $\lamstr_{v_0}\le\rho$, the following hold:
\begin{enumerate}
\item[\ttup{a}]
$\lambda_{k} <\lambda_{k-1}$ for all $k\in\NN$, unless $\lambda_k=\lamstr(\cG)$.
\item[\ttup{b}]
The sequence $\{\lambda_k\}$ converges to $\lamstr(\cG)$ as $k\to\infty$.
\item[\ttup{c}]
The sequence $\{V_k\}$ converges as $k\to\infty$, uniformly on
compact sets,
to $\Psi_*$ satisfying \cref{ET5.1A}.
\end{enumerate}
\end{theorem}

\begin{proof}
The assertions (a)--(c) follow as in the proof of
\cref{T5.3}, using \cref{L5.1,L5.2}.
\end{proof}

\begin{remark}
\cref{T5.2} describes a region of ``stability'' of the PIA using the abstract
properties in \cref{D5.2}.
We have not considered in this section  the equality
$\lamstr(\cG)=\sEmin$.
This is considered in \cref{S5.1} which follows next.
\end{remark}

%%%%%%%%%%%%%%%%%%%%%%%%%%%%%%%%%%%%%%%%%%%%%%%%%%%%%%%%%%%%%%%%%%%%%%%%%%%%%%%%
\subsection{The minimization problem for near-monotone
running costs}\label{S5.1}

We borrow the results in \cite[Proposition~5.1]{ABBK-19} which addresses
the near-monotone running cost case without imposing
any blanket stability assumptions,
and improves \cite[Proposition~1.1]{AB18}.

In general, we say that a function $f\colon\cX\to\RR$ defined on a locally compact space
is \emph{coercive, or near-monotone,
relative to a constant $\beta\in\RR$} if there exists a compact set $K$ such
that $\inf_{K^c}\,f >\beta$.

We start with the hypotheses
in \cite[Section~5]{ABBK-19} which we quote here as follows.

\begin{assumption}\label{A5.2}
In addition to the regularity hypotheses on the data
in \hyperlink{A1}{(A1)}--\hyperlink{A3}{(A3)}, we assume the following:
\begin{enumerate}
\item[(i)]
The drift $b$ and running cost $c$ satisfy, for some $\theta\in[0,1)$
and a constant $\kappa_0$, the bound
\begin{equation*}
\abs{b(x,\zeta)} \,\le\, \kappa_0\bigl(1+\abs{x}^\theta\bigr)\,,\quad\text{and\ \ }
\abs{c(x,\zeta)} \,\le\, \kappa_0\bigl(1+\abs{x}^{2\theta}\bigr)
\qquad\forall\,(x,\zeta)\in\Rd\times\Act\,.
\end{equation*}

\item[(ii)]
The drift $b$ satisfies
\begin{equation*}
\frac{1}{\abs{x}^{1-\theta}}\;
\max_{\zeta\in\Act}\;\bigl\langle b(x,\zeta),\, x\bigr\rangle^{+}
\;\xrightarrow[\abs{x}\to\infty]{}\;0\,.
\end{equation*}

\item[(iii)]
The running cost is coercive relative to $\sEmin$, that is,
\begin{equation*}
\sEmin\,<\, \lim_{r\to\infty}\,
\min_{(x,\zeta)\in\sB_r^c\times\Act} c(x,\zeta)\,.
\end{equation*}
\end{enumerate}
\end{assumption}

A full characterization of optimality under the above hypotheses
can be found in \cite[Proposition~5.1]{ABBK-19}.
We state this here in the following form.

%%%%%%%%%%%%%%%%%%%%%%%%%%%%%%%%%%%%%%%%%%%%%%%%%%%%%%%%%%%%%%%%%%%%%%%%%%%%%%%%
\begin{theorem}\label{T5.3}
Grant \cref{A5.1,A5.2}.
Then there exists a unique solution $V\in\cC^2(\Rd)$ of
\begin{equation*}
\min_{\zeta\in\Act}\, \cA^c V(x,\zeta) \,=\,
 \sEmin V(x)\quad \text{in\ } \Rd, \quad \text{and}\quad V(0)=1\,.
\end{equation*}
We also have $\sEmin=\lamstr(\cG)$, and
\ttup{a}--\ttup{c} of \cref{T3.1} hold.
 
In addition, if $c$ is near-monotone with respect to $\lamstr_v$ for
some $v\in\Usm$, then $\sE_x(c,v)=\lamstr_v$, and there exists
an inf-compact eigenfunction $\Psi_v\in\Sobl^{2,p}(\Rd)$ for $\lamstr_v$.
In particular, the diffusion controlled by such a $v$ is exponentially
ergodic.
\end{theorem}

We now state the results concerning the PIA algorithm for this model,
whose proof follows directly from \cref{T5.2,T5.3}.

%%%%%%%%%%%%%%%%%%%%%%%%%%%%%%%%%%%%%%%%%%%%%%%%%%%%%%%%%%%%%%%%%%%%%%%%%%%%%%%%
\begin{theorem}\label{T5.4}
Grant \cref{A5.1,A5.2}, suppose that
$c$ is near-monotone relative to $\rho\in(\lamstr(\cG),\infty)$, and the latter has Properties~A and B in \cref{D5.2}.
Then the conclusions of \cref{T5.2} follow.
\end{theorem}

\begin{remark}
We compare the assumptions in \cref{T5.4} to those in \cite[Theorem~5.4]{BorMey02}.
\cref{A5.2} is used to guarantee the existence of a solution
to the HJB equation, together with the standard verification of
optimality results, in \cref{T5.3}, so as to provide
a complete set of results for this model.
As remarked in that paper, \cite[Theorem~3.6]{BorMey02}
does not assert the existence of a solution to the dynamic
programming equation, but rather a dynamic programming inequality.
Existence of a solution is imposed as an assumption in the study
of the PIA.
Then, \cref{A5.1} agrees with \cite[(A4)]{BorMey02},
\hyperlink{PA}{Property~A} matches \cite[Theorem~5.4\,(i)]{BorMey02},
and \hyperlink{PB}{Property~B} is essentially the same
as \cite[Theorem~5.4\,(ii)]{BorMey02}.
\end{remark}

%%%%%%%%%%%%%%%%%%%%%%%%%%%%%%%%%%%%%%%%%%%%%%%%%%%%%%%%%%%%%%%%%%%%%%%%%%%%%%%%
%%%%%%%%%%%%%%%%%%%%%%%%%%%%%%%%%%%%%%%%%%%%%%%%%%%%%%%%%%%%%%%%%%%%%%%%%%%%%%%%
\subsection*{Acknowledgment}
We wish to thank the anonymous referee for the careful reading of
the manuscript and proposed improvements.
The research of Ari Arapostathis was supported
in part by the National Science Foundation through grant DMS-1715210, and
in part by the Army Research Office through grant W911NF-17-1-001.
The research of Anup Biswas was supported in part by
DST-SERB grants EMR/2016/004810, and MTR/2018/000028.


\end{document}